\documentclass[11pt]{amsart}

\usepackage{amsmath,amssymb,stmaryrd,dsfont}
\usepackage[pdftex]{graphicx}
\usepackage[all]{xy}

\newcommand{\var}{\varepsilon} 

\newcommand{\Vtild}{\mathop{V}\limits^{\sim\gamma}}
\newcommand{\VtildEff}{\mathop{V_{eff}}\limits^{\sim\gamma \,\,\,\,\,\,}}
 
\newcommand{\be}{\begin{equation}} 
\newcommand{\ee}{\end{equation}} 
\newcommand{\e}{\varepsilon}

\newcommand{\bac}{\begin{array}{c}}
\newcommand{\ea}{\end{array}}

\begin{document}
\newtheorem{theorem}{Theorem}[section]
\newtheorem{lemma}[theorem]{Lemma}
\newtheorem{observation}[theorem]{Observation}
\newtheorem{definition}[theorem]{Definition}
\newtheorem{example}[theorem]{Example}
\newtheorem{corollary}[theorem]{Corollary}
\newtheorem{assumption}{Assumption}
\newtheorem{property}{Property}

\numberwithin{equation}{section}

\title[Semiclassical limit for  rough Hamiltonians]{Strong and weak semiclassical limit for some rough Hamiltonians}

\author[A.  Athanassoulis]{Agissilaos  ATHANASSOULIS}
\address[A.  Athanassoulis]{DAMTP, University of Cambridge, Wilberforce Road, Cambridge CB3 0WA, UK}
\email{a.athanasoulis@damtp.cam.ac.uk}
\author[T. Paul]{Thierry Paul}
\address[T. Paul]{CNRS and CMLS, \'Ecole polytechnique, 91128 Palaiseau cedex, FRANCE}
\email{thierry.paul@math.polytechnique.fr}

%\date{}

%
%\today 
%
%
%\xxivtime
%
\begin{abstract} We present several results concerning the semiclassical
limit of the time dependent Schr\"odinger equation with potentials whose regularity doesn't
guarantee the uniqueness of the underlying classical flow.
Different topologies for the limit are considered and the situation where
two bicharacteristics can be obtained out of the same initial point is
emphasized.

\end{abstract}

\maketitle

\tableofcontents

\vskip 0.5cm

\section{Introduction}

The state of a quantum system is described by a positive, trace class operator $D^\e$, which evolves in time under the Heisenberg-von Neumann equation
\begin{equation}
\label{eq11}
\begin{array}{c}
i\varepsilon \frac{\partial }{\partial t}D^\varepsilon(t) = \left[ -\frac{ \varepsilon^2}{2}\Delta +V,
D^\varepsilon(t)\right] , \\
D^\varepsilon(t=0)=D^\varepsilon_0.
\end{array}
\end{equation}
Both the positivity and the trace of $D^\e$ are preserved in time \cite{SR}; observables, also described by operators, can be measured against the state. For example, the value of observable $\mathcal{G}^\e$, associated with operator $G^\e$, at time $t$ is given by the trace $\mathcal{G}^\e(t)=tr(G^\e D^\e(t) \,)$. The term {\em density matrix} is often used for a positive, trace-class operator.

The Heisenberg-von Neumann equation is closely related to the Schr\"odinger equation. By virtue of the singular value decomposition, the density matrix can be decomposed to a sum of {\em pure states}, i.e. rank-$1$ orthonormal projectors,
\[
D^\e=\sum\limits_m \lambda^\e_m |u^\e_m\rangle \langle\overline{u^\e_m}|
\]
 (the $\lambda^\e_m$'s are time-independent in the standard context of an external, time-independent potential, i.e. if $V$ does not depend on $D^\e$, or $t$. Here we limit our scope to such problems). Then, each of the $u^\e_m$'s satisfies a problem of the form
\begin{equation}
\label{eq1}
\begin{array}{c}
i\varepsilon \frac{\partial }{\partial t}u_m^\varepsilon 
= \left( -\frac{ \varepsilon^2}{2}\Delta +V\left( x \right) \right)
u^\varepsilon_m , \\
u^\e_m(t=0) \in L^2.
\end{array}
\end{equation}

\vskip 0.2cm

The small parameter $\varepsilon$ is called Planck's constant, and when it tends to zero one expects the system to behave in a way approaching classical mechanics; i.e. motion under Newton's law in a force field  with potential $V$. This is a singular asymptotic limit, and several techniques have been developed for its study; one is based on the Wigner transform (WT) \cite{ger,gmmp,LP}. 

The WT $W^\varepsilon(x,k,t)$ is defined as 
\begin{equation}\label{density}
 W^\varepsilon(x,k,t)=W^\e[D^\e(t)](x,k)=\int\limits_{y \in \mathbb{R}^n} { e^{-2\pi i yk} K^\varepsilon
 (x+\varepsilon \frac{y}2,x-\varepsilon \frac{y}2,t) dy }
\end{equation}
where for each $t$, $K^\varepsilon(t)$ is the integral kernel of the operator $D^\e(t)$. The WT is also closely related to the Weyl symbol of $D^\e$, and the Weyl calculus of pseudodifferential operators in general, see e.g. \cite{folland} for more details in that direction. It evolves in time under the well-known Wigner equation,
\begin{equation}
\label{eq1m2aol}
\begin{array}{l}
\partial_t {{W}^\varepsilon}(x,k) +  2\pi k \cdot \partial_x {{W}^\varepsilon}(x,k)+\\ { } \\
\,\,\,\,\,\,\,\,\,\,\,\,\,\,\,
+\frac{2}{\varepsilon} Re \left[{ 
 i \int{e^{2\pi i Sx}\widehat{V}(S){W^\varepsilon}(x,k-\frac{\varepsilon S}{2},t)dS}  }\right]=0, \\ { } \\

\,\,\,\,\,\,\,\,\,\,\,\,\,\,\,\,\,\,\,\,\,\,\,\,\,\,\,\,\,\,\,\,\,\,\,\,\,\,\,\,\,\,\,\,\,\,\,{W}^\varepsilon(t=0)={W}^\varepsilon_0,
\end{array}
\end{equation}
an equivalent reformulation of (\ref{eq11}). In this context, it can be shown that observables can be measured by
$\mathcal{G}^{\e}(t)=\int{ G_W^\e(x,k) W^\e(x,k)dxdk }$, where $G_W^\e(x,k)$ is the Weyl symbol of the operator $G^\e$. (However this in general is not a Lebesgue integral, but only makes sense as a Cauchy PV one).

In the case of a {\em pure state}, $D^\e= |u^\e\rangle \langle\overline{u^\e}|$, the WT can be written as
\begin{equation}
\label{eqyhnju}
 W^\e(x,k,t)=\int\limits_{y \in \mathbb{R}^n} { e^{-2\pi i yk}u^\varepsilon (x+\varepsilon \frac{y}2,t) 
 \bar{u}^\varepsilon (x-\varepsilon \frac{y}2,t) dy }.
\end{equation}
As a quadratic transform on an oscillating function $u^\e$, of small typical scale of oscillation $\e$, it has been widely used in signal processing, see e.g. \cite{flandrin} and a multitude of references therein for more on that direction.

Under appropriate conditions the WT  allows for a very natural and compact description of the semiclassical limit. Indeed, it has a physically meaningful  limit as $\var$ tends to zero, while in general the operator $D^\varepsilon$ itself does not. The limit (in an appropriate sense; see below for more precise statement and reference), is called the {\em Wigner measure},
\[
W^\varepsilon(t) \rightharpoonup W^0(t), 
\]
and is a probability measure that evolves in time under
the Liouville equation of classical statistical mechanics
\begin{equation}
\label{eqlolLIOUV}
\begin{array}{c}
\partial_t W^0 + 2\pi k \cdot \partial_x W^0 - \frac{1}{2\pi} \partial_x V \cdot \partial_k W^0=0, \\ { } \\
W^0_0 =\mathop{lim}\limits_{\varepsilon \shortrightarrow 0} W^\varepsilon_0.
\end{array}
\end{equation}
Observables are now described by bounded continuous functions, which can be measured against the classical state, $\mathcal{G}^0(t)=\int{G^0(x,k) dW^0(t) }$. The formal probabilistic interpretation e.g. for the expected position of a particle $\mathcal{X}(t)=\int{x dW^0(t) }$ should be intuitively clear.

Thus the correspondence between classical and quantum physics for the same system is highlighted in a very intuitive way, in terms of models that are, at least formally, similar and directly comparable. Convergence issues  are not always painless or completely understood, as we will see in more detail. There is also the well-known caveat that $W^\e(t)$ is in general neither a measure of bounded total variation, nor non-negative e.g. as a distribution. Finally asking that the limit problem itself should be well-posed isn't always trivial.

 Note at this point that considering initial data with
\[
W^\varepsilon(0) \rightharpoonup \delta(x-x_0,k-k_0)
\]
is particularly natural: they stand after all for classical particles with position $x_0$ and velocity $k_0$.

\vskip 0.2cm
It was shown in \cite{LP} (Th\'eor\`eme IV.1) that, under appropriate conditions on the initial data (which allow for practically all physically meaningful choices), as long as the potential is smooth enough, $V \in C^1$, the WT converges (in weak-$*$ sense, and up to extraction of a subsequence) to a weak solution of the limit problem (\ref{eqlolLIOUV}). It was immediately commented that, for $V \in C^1 \setminus C^{1,1}$, this does not yield an effective computation of the semiclassical limit in general, but only  a partial description: even after the extraction of a subsequence (so that there is a unique $\mathop{lim}\limits_{n} W^{\varepsilon_n}_0$), in general there are several weak solutions for the Liouville equation that $W^{\e_n}(t)$ could converge towards (because the Liouville equation with $C^1 \setminus C^{1,1}$ potential and measure valued initial data has existence, but not uniqueness of weak solutions. This can be seen e.g. using the Peano Theorem). Therefore some kind of selection principle would be needed to supplement the result, if we were to be able to identify the Wigner measure corresponding to our subsequence in $\e$ of initial data. This, in a nutshell, is the context for the present work. Our goal is to understand better the mechanisms that may lead to loss of uniqueness, and develop tools which, at least in some cases, can resolve them.

We are particularly concerned with initial data that concentrate $W^\varepsilon(0) \rightharpoonup \delta(x-x_0,k-k_0)$, and our findings here apply, under appropriate assumptions, to such data. The main assumption (other than, essentially, $C^{1, \theta}$ potential regularity), is, roughly speaking, a non-concentration condition on the flow associated with the Liouville equation.

\vskip 0.2cm
Before the statement of the present results, some reference to related works is warranted. The problem of extending semiclassical asymptotics to non-regular potentials (e.g. worse than $C^{1,1}$) has been studied by many authors -- and for good reasons: some of the simplest, most natural systems (such as a gravitational or electrostatic force field) are described by potentials that are not $C^1,1$, or $C^1$. In most cases the conclusions make use of a particular form of singularity (e.g. Coulomb/piecewise smooth potential etc) as opposed to a general smoothness class \cite{miller,ker}, and/or need some additional, non-trivial condition (non-concentration, non-interference etc) \cite{miller,FLP}. Another type of results is for whole random populations of initial data \cite{AFFGP} (which in particular can be even weaker than a result applying to ``almost all'' initial data; i.e. there might be no way to just choose  ``almost any'' initial datum at $t=0$ and keep track of it; the conclusions apply the population as a whole).  Finally, it must be mentioned that in the study of the semiclassical limit in nonlinear problems (e.g. the Schr\"odinger-Poisson equation, a very natural problem),  a lot remains to be understood. One of the main difficulties, comes precisely from the lack of regularity of the potential (which now depends on $D^\e(t)$), and the subsequent lack of well-posedness for the corresponding classical problem for measure-valued initial data; see e.g. \cite{LP,MM,ZZM}.

\vskip 0.25cm
More specifically,
semiclassical limits with rough potentials were considered recently in 
\cite{AFFGP} where the hypotheses on the potential are, roughly speaking, that it has a generic part with $BV_{loc}$ gradient, plus possibly a {\em repulsive} Coulomb part.
The setting involves a random population of initial data, and it is shown that the corresponding population of solutions at a later time   tends weakly to the push-forward by the Liouville equation of a population of measures. The classical propagation problem for the random family of measures is shown to be well-posed by virtue of the Ambrosio-Di Perna-Lions theory \cite{ambrosio,ldp}. In \cite{FLP} similar core ideas were used with deterministic mixed states. In
that context, an averaging condition forbids  initial 
%(density matrix) 
data concentrating e.g. to a delta function in phase-space, a fact consistent with
the weak limit strategy and the fact that the flow is not defined everywhere.
The present work therefore concerns less general  potentials than those treated in \cite{AFFGP,FLP}, (still giving rise to ill-posed classical dynamics of course), but can deal with initial data concentrating 
to point-supported measures.

The possible ill-posedness of the limit equation can originate not only from lack of smoothness of the potential, which we focused on for our discussion so far, but e.g. from working on a torus instead of Euclidean space \cite{mac}.

\vskip 0.5cm
Finally, let us present a concrete example.
Already in \cite{LP} ,
families of initial data were constructed so that their Wigner measures at later times oscillate between the several solutions of the
Liouville equation. 
The one-dimensional case study used for that, was essentially equivalent to the following:

Consider the potential
\be\label{vlp}
V(x)=-\vert x\vert^{1+\theta}\,\cdot\,\beta(x)
\ee
with $\theta\in(0,1),\ C^\infty\ni\beta=1$ on $[-1,+1]$, and $\beta(x)=0$ for $\vert x\vert\geq 2$, and initial data such that $W^\e_0 \rightharpoonup \delta(x-0,k-0)$. Then one can
easily check that there exist, for $t$ small enough, two families of bicharateristics starting near the origin at $t=0$, namely
\begin{equation}\label{bichar}
(X^\pm(t), P^\pm(t)) = ( \pm c_0 t^\nu ,\pm\frac{c_0\nu}{2\pi} t^{\nu-1}  ), 
\ee
with  $\nu=\frac{2}{1-\theta}$ and $c_0=\left({  \frac{(1-\theta)^2}{2} }\right)^{1-\theta}$.

In this case it is possible to find semiclassical limits where one wavepacket (a $\delta$-function in classical phase-space) splits in two wavepackets. (In this particular example it happens immediately, i.e. $\forall t>0$, or not at all; different examples can be constructed where it happens after an $\e$-independent time $t_*$). The mass of each wavepacket depends on information which can be easily extracted from the quantum problem, but is lost if we take the limit as usual in a straightforward way.

\vskip 0.5cm

The first of our results deals explicitly with the  case study from \cite{LP} mentioned earlier, to provide a clear illustration to the selection principle at work. 
Our second Theorem gives a general result concerning the approximation, in strong
topology, of the Wigner function of the solution of (\ref{eq1m2aol}) with the solution of a ``smoothed'' Liouville equation. 
Our third results goes back to the weak approximation, but with explicit remainder estimates.

\section{Main results}
\label{subseczoro}

We will use  a  space of test-functions   introduced in \cite{LP}: let
\begin{equation}
\begin{array}{c}
\label{eqzato}
\mathcal{A} = \{ \phi \in C_0(\mathbb{R}^{2n}) \,|\,\, \int{ \mathop{sup}\limits_{x} |\mathcal{F}_{k \shortrightarrow K}[\phi(x,k)]|dK } < \infty \}, 
\end{array}
\end{equation}
equipped with the norm
\begin{equation}
\begin{array}{c}
||\phi||_{\mathcal{A}}= \int{ \mathop{sup}\limits_{x} |\mathcal{F}_{k \shortrightarrow K}[\phi(x,k)]|dK }. 
\end{array}
\end{equation}

\vskip 0.35cm
Our first result will focus on the case study presented in the introduction; a  family of initial data which concentrates at
the origin at time $t=0$ and splits into two separated wave packets at time $t>0$. This is only an example among several
others
presented in section \ref{secExamples} and the more general results of the paper are given by the Theorems 
\ref{DaTheo} and \ref{THEO22} below.

\begin{theorem}\label{thetheo}
Let $V$ given by \eqref{vlp} and
\be\label{ini}
W_0^\varepsilon(x,k)=\lambda^{\frac{7+3\theta}{30}}w( {\lambda^{\frac{1+\theta}{6}}}x,
 {\lambda^{\frac{1-\theta}{15}}}k) \ast \left({\frac{2}{\var}}\right)^n e^{-2\pi\frac{ x^2+k^2 }{\var}}
\ee
with $\lambda=log(\frac{1}\e)$, $w \in H^2 \cap L^\infty \cap L^1$, $w(x,k)\geqslant 0$, 
 $supp\, w \subseteq \{ |x|^2+|k|^ 2<1 \}$, $\int{w(x,k)dxdk}=1$.

Then $\exists T>0$ such that for all $t\in[0,T]$, the solution $W^\e(t)$ of \eqref{eq1m2aol}
converges in weak-$*$ sense in $\mathcal{A}'$ to
\be\label{split}
W^0(t)=c_+\delta{(X^+(t),P^+(t))}+c_-\delta{(X^-(t),P^-(t))},
\ee
with $(X^\pm(t),P^\pm(t))$ given by \eqref{bichar} and
\[
c_\pm = \int\limits_{\pm x >0}{w(x,k)dxdk}.
\]
\end{theorem}

\noindent {\bf Remark:} The fact that the initial data contain a smoothing by $\left({\frac{2}{\var}}\right)^n e^{-2\pi\frac{ x^2+k^2 }{\var}}$ ensures that $W_0^\var$ corresponds to a density matrix, i.e. a positive trace-class operator (see e.g. Exemple III.7 of \cite{LP}, or Lemma \ref{dktnkmmtl}).

\vskip 0.4cm

%Let us compare this result with the results in \cite{LP} concerning the potential \eqref{vlp}. (There are  differences in the setup, but for small enough time, i.e. until the characteristics leave $\{ |(x,k)|<1 \}$, the two results are directly comparable). 

In Remarque IV.3 of \cite{LP}, one can find a construction of a family of normalized wave functions $u_\e^+$ concentrating at the origin, and whose
Wigner functions, after extraction of a subsequence, converge to a Dirac mass centered on $(X^+(t), P^+(t))$. The same construction is obviously possible in order to get a family $u_\e^-$ whose
Wigner functions, after extraction of a subsequence, converge to a Dirac mass centered on $(X^-(t), P^-(t))$.
Defining $D^\e_0=\frac1{\sqrt2}(\vert u_\e^+\rangle\langle u_\e^+\vert+\vert u_\e^-\rangle\langle u_\e^-\vert)$ 
one can  check that the Wigner function of $D^\e_0$ will follow the conclusion of Theorem \ref{thetheo}.  
The present construction is different, and is not obtained through compactness arguments: there is no need of extraction of subsequence, and the scaling property is explicit (it is given by an implicit diagonal argument in a two scale sequence in \cite{LP}). 
%Moreover, as shown in section \ref{secExamples} the present construction extends to more general examples.
The construction of \cite{LP} is based on the more-or-less explicit understanding of the flow around the singular point, while our approximation (i.e. Theorems \ref{DaTheo} and \ref{THEO22} below) is built around a 
non-concentration condition, which is relatively easy to check for other problems, whether the flow can be computed explicitly or not. (This may be increasingly important in high dimensions).

Let us also point out the similarity with the long-time behaviour around a regular separatrix, studied extensively in \cite{tp2}. In the regular case, the quantum wavepacket leaves the fixed point in ``infinite'' (i.e. asymptotically long) time, as does the classical flow; here the classical flow leaves the fixed point immediately, and this behaviour (under additional assumptions, of course) is preserved by the quantum system.

\vskip 0.5cm

We turn now to the more general results of this paper.
\vskip 0.5cm

In the sequel we will consider potentials $V$  that satisfy the following assumption:
\begin{assumption} \label{cond1}%\nonumber 
\begin{equation}
\int_{\mathbb{R}^{n}}|\widehat{V}(S)| \,\, \frac{S^2}{1+S^2} \,\,dS < \infty,
\end{equation}
 and moreover 
there are constants $C>0, \theta \in (0,1)$ such that for $m \in \{0,1,2\}$
\[
\begin{array}{c}
\forall 1\leqslant a\leqslant b \leqslant +\infty : \int\limits_{|S|\in (a,b)}{|\widehat{V}(S)| \, |S|^m dS} 
\leqslant \frac C{m-1-\theta}\left(b^{m-1-\theta}-a^{m-1-\theta}\right)
%\int\limits_{\rho \in (a,b)}{ \rho^{-2-\theta+m}d\rho} = \\ { } \\
%= \left.{\frac{\rho{-1-\theta+m}}{-1-\theta+m}}\right]_{a}^{b}
\end{array}
\]
\end{assumption}
This is closely modeled after $V(x)=C |x|^{1+\theta}$; indeed it is easy to check that the aforementioned potential satisfies this condition. In section \ref{secExamples}, we will also see  some other relevant types of singularities (also generated from $|x|^{1+\theta}$ in some sense).

\vskip 0.25cm

Denote
\begin{equation} \label{eqdefVtilde}
\Vtild(x) = \left({ \frac{2}{\e^{\gamma } } }\right)^{\frac{n}{2}} \int{ e^{-\frac{2\pi}{\e^\gamma  }|x-x'|^2} V(x')dx' }.
\end{equation}

\vskip 0.2cm

\begin{theorem}\label{DaTheo}
Let us suppose that Assumption \ref{cond1} holds, and there exist $T>0$, $\delta \in (0,\frac{\theta}{2+\theta})$, $\gamma>\frac{2}{1+\theta}$ and approximate initial data $[W_0^\varepsilon] \in H^2(\mathbb{R}^{2n})$ such that 
\begin{equation}\label{approxi}
||[W_0^\varepsilon]-W_0^\varepsilon||_{L^2} =o(||W_0^\varepsilon||_{L^2})
\end{equation}
and the solution of
\begin{equation}
\label{eqmaikol22}
 \partial_t \rho+2\pi k \cdot \partial_x \rho-\frac{1}{2\pi} \partial_x \Vtild \cdot \partial_k \rho =0,
\end{equation}
with initial condition $\rho(t=0)=[W_0^\varepsilon]$
satisfies 
\begin{equation}
\label{paoeqmmiikkooll}
||\rho(t)||_{H^2} =O(\varepsilon^{-\delta}||W_0^\varepsilon||_{L^2})
\end{equation}
uniformly on $[0,T]$.
  \vskip 0.5cm
  
Then the Wigner function $W^\e(t)$ of the solution of  (\ref{eq1m2aol}) satisfies, uniformly on $[0,T]$,
%\eqref{eq1} 
\begin{equation}
\label{polai8d}
|| W^\e(t)-\rho_1^\varepsilon(t)||_{L^2} =O(\varepsilon^{\kappa}||W_0^\varepsilon||_{L^2} \,\, + \,\, ||W_0^\varepsilon-[W_0^\varepsilon]||_{L^2})=o(||W^\e_0||_{L^2}).
\end{equation}
where $\rho_1^\e$ is the solution of \eqref{eqmaikol22} with initial datum $\rho_1^\varepsilon(t=0)=W_0^\varepsilon$ and
$\kappa=min\{ \,\, \gamma\frac{1+\theta}{2}-1, \,\, \frac{\theta}{2+\theta}-\delta \,\,\}$.
\end{theorem}

\vskip 0.2cm
\noindent {\bf Remark: The non-concentration condition.} It is helpful to discuss the motivation behind this type of assumption. Ideally, one would like the solution of the problem
\[
\bac
 \partial_t \rho+2\pi k \cdot \partial_x \rho-\frac{1}{2\pi} \partial_x \Vtild \cdot \partial_k \rho =0, \\
 \rho(t=0)=W_0^\e
\ea
\]
to have $H^2$ regularity in the sense $||\rho(t)||_{H^2}\leqslant e^{Ct} ||W_0^\e||_{H^2}$ for some $C$ that is either $\e$-independent, or grows extremely slowly in $\e$. Of course, when working with low-regularity potentials, this is in general not possible. So one looks for a weaker version of regularity, which would still be sufficient here; it turns out that the property
\[
\bac
\exists [W_0^\varepsilon] \in H^2(\mathbb{R}^{2n}), \,\,\,\,\,\,\,\, ||[W_0^\varepsilon]-W_0^\varepsilon||_{L^2} =o(||W_0^\varepsilon||_{L^2}) \,\,\,\, \mbox{such that, if} \\
 \partial_t \rho+2\pi k\cdot \partial_x \rho-\frac{1}{2\pi} \partial_x \Vtild \cdot \partial_k \rho =0,  \,\,\,  \rho(t=0)=[W_0^\e] \\
\mbox{then }  \,\, ||\rho(t)||_{H^2} \,\, \mbox{ is not too large}
\ea
\]
is an acceptable surrogate. 

We call this a non-concentration condition, because of one of the simplest ways to check it. Assume that there is a ``small'' set  $S \subseteq \mathbb{R}^{2n}$, such that the derivatives of $\Vtild$ up to order $3$ outside $S$ are ``not too large'' in $\e$. Now set $B=\mathop{\bigcup}\limits_{t \in [0,T]} \phi^\e_{-t} (S)$, where $\phi^\e_t$ is the flow associated with the regularized Liouville equation (\ref{eqmaikol22}). $B$ contains all the ``dangerous'' part of phase-space; away from $B$ the original regularity requirement holds automatically. A ``concentrating flow'' could pull too big a part of phase-space through the original ``small bad neighbourhood'' $S$ in time $t \in [0,T]$. On the other hand, a ``repulsive singularity'', in general tends to push the flow away, thus keeping the pre-image of $S$ small enough. The approximate initial data $[W_0^\e]$ can be then constructed by restricting the initial data on $B^C$.

\vskip 0.2cm 
\noindent {\bf Remark: Checking the assumptions for concrete problems.} Concrete examples for which the assumptions (in particular the non-concentration condition) can be checked to hold, are examined in section \ref{secExamples}. The techniques used can be easily extended to similar problems. It must be noted that the condition for Theorem \ref{DaTheo} is seen to be satisfied for a generic selection of initial data, as long as the rate of concentration satisfies certain constraints (i.e. the profile that concentrates can be any function with a given regularity; see e.g. Lemmata \ref{lmscalingsonso}, \ref{lmtbswpolmlkbl}). In contrast, for the stronger Theorem \ref{THEO22} (see next Remark for the difference), in general we need to restrict the shape of the wavepacket with additional assumptions for the same problems; Lemmata \ref{lpbuse9876v}, \ref{lemthirdartbg}.

\vskip 0.2cm
\noindent {\bf Remark: $L^2$ versus measures, and convergence.} Theorem \ref{DaTheo} is an $L^2$ asymptotic approximation result: two functions, of constant $L^2$ norm in time (also equal to one another and large in $\e$), are relatively close in $L^2$. It is not, however,  a convergence result, since in general neither function has a limit in $L^2$. Ideally, we would like to conclude that, since they're (relatively) close, both functions converge (in some appropriate sense) to the same limit, typically a (sum of) delta function(s). However,  the (relative) $L^2$ approximation is not by itself  sufficient to ensure that if e.g. $\rho_1^e(t) \rightharpoonup \delta(x-x(t),k-k(t))$, then also $W^\e(t)$ goes to the same limit. Is there some additional information in the problem that we could use to come to that conclusion? Or is there a physical issue here, of $L^2$ norm really ``missing'' the solution? For example, is it really possible that in a problem with splitting (as in Theorem \ref{thetheo}) the conclusion of Theorem \ref{DaTheo} holds,  in the limit $\rho^\e_1$ ``goes to the right'', but at the same time $W^\e$ (or half its mass) ``goes to the left''? This is a question that we cannot really answer at the moment. One possible way to go is by strengthening the non-concentration condition, and that is pursued in Theorem \ref{THEO22} below. This yields full control over the convergence (in an appropriate sense, see also Corollary \ref{coro1}), but at the cost of excluding many interesting, natural problems.

\vskip 0.2cm

\begin{theorem} \label{THEO22} 
Assume that the hypotheses of Theorem \ref{DaTheo} hold and $\rho_1^\varepsilon(t)$ is as in Theorem \ref{DaTheo}.  
Moreover, in addition to the condition (\ref{approxi}) we assume that
\begin{equation}\label{approxi2}
||[W_0^\varepsilon]-W_0^\varepsilon||_{L^2} =o(1),
\end{equation}
$||W_0^\var||_{L^2}=o(\var^{-\kappa})$ ($\kappa$ being the same as in Theorem \ref{DaTheo}), and the initial data $W_0^\var$ corresponds to a density matrix, i.e. it is the Wigner function of a positive trace-class operator $D_0^\var$, $tr(D^\var_0)=1$ $\forall \var>0$.

\vskip 0.25cm
For $\phi \in  \mathcal{A}$ denote
\[
 E^\varepsilon(\phi)
= \langle W^\varepsilon(t)-\rho_1^\varepsilon(t), \phi \rangle .
\]
and, for brevity, by $F(\var)$ the bound of equation (\ref{polai8d}), namely
\[
F(\var)=C( \var^\kappa||W_0^\var||_{L^2}+||[W_0^\varepsilon]-W_0^\varepsilon||_{L^2})=o(1),
\]
for some large enough constant $C>0$.

\vskip 0.25cm
Then, $\forall \phi \in  \mathcal{A}, \, t\geqslant 0,$ the following estimate holds:
\begin{equation}
\label{eq214}
\begin{array}{c}
|E^\varepsilon(\phi)| \leqslant F(\var)||\phi||_{\mathcal{A}} + \frac{2^n \pi^{\frac{n}2}}{\sqrt{\Gamma(n+1)}} \sqrt{F(\var)} +\\ { } \\

 + (1+||[W_0^\var]||_{L^1}) \,\,\,\,  ||\,(1-\chi_{F(\var)^{\frac{1}{2n}}}) \, \phi \,||_{\mathcal{A}}

\end{array}
\end{equation} 
where $\chi\in C^\infty(\mathbb{R}^{2n},[0,1])$ is any cutoff function, $\chi(z) =1 \,\, \forall |z|<1$, $\chi(z) =0 \,\, \forall |z|>2$, and $\chi_M(z):=\chi(Mz)$.

\end{theorem}

\vskip 0.2cm
\noindent {\bf Remark: $L^2$ versus measures, and convergence.} If a generic ``nice'' function concentrates, $f^\e \rightharpoonup \delta(x,k)$, and $||f^\e-[f^\e]||_{L^1}=o(1)$, then it follows that $[f^\e] \rightharpoonup \delta(x,k)$ as well. In other words $||[W_0^\varepsilon]-W_0^\varepsilon||_{L^1} =o(1)$ (which in general holds in the examples where Theorem \ref{DaTheo} does) would better suited to control concentration. We are forced however to use the much stronger $||[W_0^\varepsilon]-W_0^\varepsilon||_{L^2} =o(1)$ due to technical reasons, and in particular the fact that the Wigner equation doesn't respect the $L^1$ norm in general (or the total variation of measures, or even the $\mathcal{A}'$ norm away from the positive cone which also scale like $L^1$). The assumptions of Theorem \ref{THEO22} should be considered pessimistic, and  probably can be relaxed. In any case, the precise convergence result is summarized in Corollary \ref{coro1} below.

\vskip 0.2cm
\noindent {\bf Remark: Size of the error.} The last term in rhs of equation (\ref{eq214}) can be further estimated if there is more information on $\phi$, and in particular its decay in the $k$ variable.
\vskip 0.25cm

\begin{corollary}\label{coro1} If the assumptions of Theorem \ref{THEO22} are satisfied and $||[W_0^\var]||_{L^1}$ is uniformly bounded in $\var$, then
\[
\mathop{lim}\limits_{\var \shortrightarrow 0} E^\var(\phi)=0 \,\,\,\,\,\,\, \forall \phi \in  \mathcal{A}, \, t\geqslant 0,
\]
and therefore the Wigner function $W^\varepsilon(t)$ has a weak-$*$ limit in $\mathcal{A}'$ 
 if and only if $\rho_1^\varepsilon$ has a weak-$*$ limit, and the two are equal, i.e.
\begin{equation}
\exists W^0(t) \,\,:\,\,  \mathop{lim}\limits_{\varepsilon \shortrightarrow 0} \langle W^\varepsilon(t), \phi \rangle
=\langle W^0(t), \phi \rangle \,\,\, \forall \phi\in\mathcal{A}, \,\,  t \in [0,T],
\end{equation}
if and only if
\begin{equation}
 \mathop{lim}\limits_{\varepsilon \shortrightarrow 0} \langle \rho_1^\varepsilon(t), \phi \rangle=\langle W^0(t), \phi \rangle \,\,\, \forall \phi\in\mathcal{A}, \,\,  t \in [0,T].
\end{equation}
\end{corollary}

In other words, under the assumptions of Theorem \ref{THEO22}, one has to solve equation (\ref{eqmaikol22}) in order to compute the propagation in time, $\rho_1^\e(t)$, and then take the limit $\e \shortrightarrow 0$.  That is, unlike what happens with regular potentials, the semiclassical limit and the propagation in time by (a regularization of) the Liouville flow no longer commute.

\begin{figure}[htb!]
\hskip 0.1cm
\xymatrix{
\rho_1^\varepsilon(0)  \ar[r]^{\circ \phi^\varepsilon_t} \ar[d]_{\varepsilon \shortrightarrow 0} & \rho_1^\varepsilon(t) \ar[d]^{\varepsilon \shortrightarrow 0}   \\
 W^0(0)  & W^0(t) } 
\caption{In some cases, there is not enough information in $W^0(0)$ to determine $W^0(t)$, i.e. the Cauchy problem (\ref{eqlolLIOUV}) is ill-posed.}
\label{Fig200}
\end{figure}

\vskip 0.5cm

\noindent {\bf Organization of the paper:}  Proofs of the main Theorems are in section \ref{987koiu}. Checking the assumptions of the Theorems for various concrete problems is done in section \ref{secExamples}. Certain auxiliary results we use are collected in section \ref{secauxres}.

\section{Definitions and Notations}
\label{secDef}

The Fourier transform is defined as
\begin{equation}
\widehat{f}(k)=\mathcal{F}_{x \shortrightarrow k} \left[ {f(x)} \right]
=\int\limits_{x\in \mathbb{R}^n} {e^{-2\pi i k x} f(x)dx}. 
\end{equation}
%Inversion is given by
%\begin{eqnarray}
%\check{f}(k)=\mathcal{F}^{-1}_{x \shortrightarrow k}\left[f(x)\right]
%=\int\limits_{x\in \mathbb{R}^n} {e^{2\pi ikx}f(x)dx}, \\
%\mathcal{F}_{b \shortrightarrow x}\left[ \mathcal{F}^{-1}_{a \shortrightarrow b}\left[ f(a) \right] \right]=
%\mathcal{F}^{-1}_{b \shortrightarrow x}\left[ \mathcal{F}_{a \shortrightarrow b}\left[ f(a) \right] \right]=
%f(x).
%\end{eqnarray}
%
\vskip 0.2cm

 For compactness, we will use the following notations:
\[
\begin{array}{c}
T^V_\var W = \frac{2}{\varepsilon} Re \left[{ 
 i \int{e^{2\pi i Sx}\widehat{V}(S){W}(x,k-\frac{\varepsilon S}{2})dS}  }\right]=
\\ { } \\

=2\mathcal{F}^{-1}_{X,K \shortrightarrow x,k} \left[{ \int{\widehat{V}(S)\widehat{W}(X-S,K)\frac{sin(\pi \var SK)}{\var}dS} }\right],
\end{array}
\]
\vskip 0.1cm
\[
\begin{array}{c}
T^V_0 W=-\frac{1}{2\pi} \partial_x \Vtild \cdot \partial_x W=\\ { } \\

=2\pi \mathcal{F}^{-1}_{X,K \shortrightarrow x,k} \left[{  \int{\widehat{V}(S)\widehat{W}(X-S,K)S \cdot K dS} }\right].
\end{array}
\]

\vskip 0.2cm

The Sobolev norms of order $m$ on phase-space will be defined as follows:
\[
||f||_{W^{m,p}(\mathbb{R}^{2n})}= \sum\limits_{|a|+|b| \leqslant m} ||\partial_x^a \partial_k^b f ||_{L^p(\mathbb{R}^{2n})},
\]
where of course $a$ and $b$ are multi-indices of length $n$ each. Moreover, $H^m(\mathbb{R}^{2n})=W^{m,2}(\mathbb{R}^{2n})$.

\section{Proofs} \label{987koiu}

\subsection{Proof of Theorem \ref{thetheo}} \label{sspoftlft} 
In this section we will show that Theorem \ref{thetheo} follows from Theorems \ref{DaTheo}, \ref{THEO22}. The latter are proved in the next subsections.

\vskip 0.1cm
Before getting into more detail, let us comment that $\lambda$ as used in the statement of Theorem \ref{thetheo} is an effective upper bound. It is easily checked that if we replace $\lambda$ by  $\lambda' \leqslant \lambda,\ \frac 1{\lambda'}=o(1)$, the Theorem still holds.

\noindent {\bf Remark on notation:}  {\em
In this section we use the notation $R=\lambda^{-\frac{1}{3}}$.  }

\vskip 0.2cm
Now denote
\be\label{inisss}
f_0^\varepsilon(x,k)=\lambda^{\frac{7+3\theta}{30}}w( {\lambda^{\frac{1+\theta}{6}}}x,
 {\lambda^{\frac{1-\theta}{15}}}k),
\ee
where of course $\lambda=log(\frac{1}\e)$,
i.e. $W^\var_0=f^\var_0   \ast \left({\frac{2}{\var}}\right)^n e^{-2\pi\frac{ x^2+k^2 }{\var}}$. One should note that 
\be
||W^0_\var-f_0^\var||_{L^2}\leqslant \sqrt{\var} ||f^\var_0||_{H^1} =  O(\var^{\frac{1}2-\eta}) \,\,\,\, \forall \eta\in (0,\frac{1}{2})
\ee
(The approximation follows from a standard observation on smoothing operators, see e.g. \cite{AP} for a proof.
Moreover, since by construction the constant $\kappa$ of Theorem \ref{DaTheo} is $\kappa < \frac{1}3$, the error $||W^0_\var-f_0^\var||_{L^2}$ can indeed be dropped without loss of generality, when invoking Theorem \ref{THEO22} -- and Theorem \ref{DaTheo} through it).

%$\lambda$ as stated is an effective upper bound; any $1 \ll \lambda' \leqslant \lambda$ is also admissible. 
%For the computation of the constant in the definition of $\lambda$ (which is of course very crude, but nevertheless explicit, so as not to clutter the statement with unimportant auxiliary parameters) see observation \ref{hjancpalmkl}.

\vskip 0.3cm
In Lemma \ref{lpbuse9876v} it is shown essentially that Theorem \ref{THEO22} applies if $n=1$, $V(x)=-|x|^{1+\theta}$ in $\{ |x|<1 \}$ (with a smooth cutoff outside of that, the specifics of which are irrelevant) and $W_0^\varepsilon=\delta_x ^{-1} \delta_k^{-1} w(\frac{x}{\delta _x},\frac{k}{\delta_k})$, where $\delta_k = R^\frac{1+\theta}{2}$, $\delta_x = R^{\frac{1-\theta}{5}}$. (We avoid duplicating parts of the proof of Lemma \ref{lpbuse9876v} here. Moreover, it is completely straightforward to see that the assumptions not checked explicitly there also hold).  Therefore Corollary \ref{coro1} applies, i.e. it suffices to find the limit in $\e$ of $\rho_1^\e$.

\vskip 0.2cm
\noindent {\bf Claim:} With $[W_0^\e]$ defined as in equation (\ref{eqjkjktndh}), we have
\begin{equation}
\begin{array}{c}
||W_0^\varepsilon-[W_0^\varepsilon]||_{L^1} =o(1).
\end{array}
\end{equation}

\vskip 0.2cm
\noindent {\bf Proof of the claim:} We will work separately for the two terms
\[
||W_0^\varepsilon-[W_0^\varepsilon]||_{L^1} \leqslant ||W_0^\varepsilon-f_0^\varepsilon||_{L^1} +||f_0^\varepsilon-[W_0^\varepsilon]||_{L^1}.
\]
First of all, it is obvious that 
\[
\int\limits_{|(x,k)|>2} { |f^\var_0-W_0^\var|dxdk }=\int\limits_{|(x,k)|>2} { |W_0^\var|dxdk }=O(\var^\infty).
\]
Moreover, using the Jensen inequality, one sees that
\[
\begin{array}{c}
\int\limits_{|(x,k)|<2} { |f^\var_0-W_0^\var|dxdk } \leqslant \sqrt{ \int\limits_{|(x,k)|<2} {dxdk} }
\sqrt{\int\limits_{|(x,k)|<2} { |f^\var_0-W_0^\var|^2dxdk }} \leqslant \\ { } \\

\leqslant C||f^\var_0-W_0^\var||_{L^2},
\end{array}
\]
and therefore
\[
||W_0^\varepsilon-f_0^\varepsilon||_{L^1} \leqslant O(\var^\infty) + C||f^\var_0-W_0^\var||_{L^2}=o(1).
\]
Moreover,
\begin{equation}
\begin{array}{c}
||f_0^\varepsilon-[W_0^\varepsilon]||_{L^1} \leqslant C \delta_x^{-1} \delta_k^{-1} \int\limits_{|x|<CR} {|w(\frac{x}{\delta_x},\frac{k}{\delta_k})|dxdk}= \\ { } \\

=  \int\limits_{|x|<C\frac{R}{\delta_x}} {|w(x,k)|dxdk} \leqslant C
\delta_x^{-1} R =o(1).
\end{array}
\end{equation}
The proof of the claim is complete.

This means that $\rho_2^\e$ and $\rho_1^\e$ (the evolution under equation (\ref{eqmaikol22}) of $[W_0^\e]$, $W_0^\e$ respectively) are interchangeable for our purposes, since
\[
\langle \rho_1^\e,\phi \rangle = \langle \rho_2^\e,\phi \rangle + \langle \rho_1^\e-\rho_2^\e,\phi \rangle
\]
and
\[
|\langle \rho_1^\e-\rho_2^\e,\phi \rangle| \leqslant ||\rho_1^\e-\rho_2^\e||_{L^1} ||\phi||_{L^\infty} \leqslant  || W^\e_0-[W^\e_0] ||_{L^1} ||\phi||_{\mathcal{A}}.
\]
(We used the obvious bound $||\phi||_{L^\infty} \leqslant  ||\phi||_{\mathcal{A}}$; recall that $\mathcal{A}$ was defined in eq. (\ref{eqzato}) ).

\vskip 0.5cm
In particular, this means we can always work with $|x|>CR$, since $\rho_2^\e$ stays outside of $\{ |x|<CR\}$ by construction. This allows certain ODEs we will use to be well-posed. 
Without loss of generality we will work for $x>0,k>0$ (working with $k<0$ makes no difference other than the opposite sign in the explicit expression for $K(0)$ in equation (\ref{piupiuequation}), 
the case $x<0$ follows by symmetry). Denote by $X(t),K(t)$ the unique solution of 
\begin{equation}
\begin{array}{c}
\label{piupiuequation}
\dot{X}(t)=2\pi K(t), \,\,\,\,\,\,\,\,
\dot{K}(t)=\frac{1+\theta}{2\pi} \left({X(t)}\right)^\theta, \\ { } \\

X(0)=C R, \,\,\,\, K(0)=|X(0)|^{\frac{1+\theta}{2}},
\end{array}
\end{equation}
i.e. the only branch of the level set of $\{ 2\pi^2k^2+V(x)\}=0$ in $\{ x>R \,\,\wedge \,\, k>0 \}$. 

 The interesting property of this trajectory is that, unlike what happens with regular potentials, it leaves zero in finite time (more generally, i.e. without restricting to $\{ x>R \,\,\wedge \,\, k>0 \}$, it reaches {\em and} leaves zero in finite time). 

\noindent {\bf Claim:} {\em Any characteristic of the Liouville equation (\ref{eqmaikol22}) starting in $\{ x>R \,\,\wedge \,\, k>0 \} \cap supp [W_0^\e]$ converges to $X(t),K(t)$ defined in equation (\ref{piupiuequation}).} 

Then it readily follows that 
\[
\chi_{x>0}(x,k)[W_0^\e] \rightharpoonup \int{\chi_{x>0}(x,k)w(x,k) dxdk} \,\, \delta(x-X(t),k-K(t)).
\]
 Repeating the argument for $x<0$ gives a limit of two delta functions leaving zero along different trajectories.

\vskip 0.2cm

\noindent {\bf Proof of the claim:} Denote by $X_1(t),K_1(t)$ the solution of
\begin{equation}
\begin{array}{c}
\label{piupiuequation2}
\dot{X}_1(t)=2\pi K_1(t), \,\,\,\,\,\,\,\,
\dot{K}_1(t)=-\frac{1}{2\pi} \partial_x \Vtild(X_1(t)), \\ { } \\

X_1(0)=x_0, \,\,\,\, K_1(0)=k_0,
\end{array}
\end{equation}
where of course $(x_0,k_0) \in \{ x>R \,\,\wedge \,\, k>0 \} \cap supp [W_0^\e]$ are as in the statement of the claim). Then equation (\ref{piupiuequation2}) can be recast as
\begin{equation}
\begin{array}{c}
\label{piupiuequation3}
\dot{X}_1(t)=2\pi K_1(t), \,\,\,\,\,\,\,\,
\dot{K}_1(t)=-\frac{1}{2\pi} {V}_x(X_1(t)) - \frac{1}{2\pi}(\partial_x\Vtild(X_1(t))-{V}_x(X_1(t))), \\ { } \\

X_1(0)=x_0, \,\,\,\, K_1(0)=k_0.
\end{array}
\end{equation}
Here we use the estimate $|\partial_x\Vtild(X_1(t))-{V}_x(X_1(t))| \leqslant C R^{\theta-1} \e^{\frac{\gamma}{2}}=o(1)$ (see Lemma \ref{auxlemtkkors} with $V'(x)$ in the place of $f$). It is clear now that any such characteristic converges to
\begin{equation}
\begin{array}{c}
\label{piupiuequation4}
\dot{X}_2(t)=2\pi K_2(t), \,\,\,\,\,\,\,\,
\dot{K}_2(t)=-\frac{1}{2\pi} {V}_x(X_2(t)) , \\ { } \\

X_2(0)=x_0, \,\,\,\, K_2(0)=k_0.
\end{array}
\end{equation}
To conclude observe that trajectory $X_2(t),K_2(t)$ is squeezed between two level sets of $2\pi^2k^2+V(x)$ converging to each other (namely between $\{ 2\pi^2k^2+V(x)\}=V(CR)$ and $\{ 2\pi^2k^2+V(x)\}=V(R^{\frac{1-\theta}5})$). 

\vskip 0.5cm
Essentially the same analysis applies to the example of Lemma \ref{lemthirdartbg} as well.

\subsection{Proof of Theorem \ref{DaTheo}} 
\label{secProofs}

\vskip 0.2cm

It will be helpful to recall the main objects we are going to use here. The WT for this problem, $W=W^\var(x,k,t)$, satisfies the well-known Wigner equation (\ref{eq1m2aol}) with initial data $W^\e_0$.
%\begin{equation}
%\label{WTeq}
%\left\{ {
%\begin{array}{c}
%\partial_t {W^\e} +  2\pi k \cdot \partial_x {W^\e}+ T^V_\var W^\e=0, \\ { } \\
%
%W^\e(t=0)={W}^\varepsilon_0.
%\end{array}
%}\right.
%\end{equation}
Moreover, $W^\e_1$ and $\rho_2^\e$ are defined as follows
\begin{equation}
\label{W1Eq}
\left\{ {
\begin{array}{c}
\partial_t {W^\e_1} +  2\pi k \cdot \partial_x {W^\e_1}+ T^{\Vtild}_\var W^\e_1=0, \\ { } \\

{W^\e_1}(t=0)={W}^\varepsilon_0,
\end{array}
}\right.
\end{equation}
and $\rho_2^\e$,
\begin{equation}
\label{eqmaikol1922}
\left\{ {
\begin{array}{c}
 \partial_t \rho_2^\e+2\pi k \partial_x \rho_2^\e+T^{\Vtild}_0 \rho_2^\e =0,\\ { } \\
\rho_2^\e(t=0)=[W_0^\varepsilon],
\end{array}
}\right.
\end{equation}
where $[W_0^\varepsilon]$ for now is simply assumed to exist (and have the properties found in the statement of Theorem \ref{DaTheo}. Such approximate initial data will be constructed, for concrete examples, in Section \ref{secExamples}).
Finally, recall that $\rho_1^\varepsilon$ was defined as the solution of equation (\ref{eqmaikol22}) with initial data $\rho_1\e(t=0)=W^\e_0$.

\vskip 0.25cm
We partition the proof as follows:
\begin{equation}
\begin{array}{l}
||W^\e-\rho_1^\varepsilon||_{L^2} \leqslant ||W^\e-W^\e_1||_{L^2} + \\ { } \\

\hskip 1cm +||W^\e_1-\rho_2^\e||_{L^2}+||\rho_1^\varepsilon-\rho_2^\e||_{L^2} =o(||W_0^\varepsilon||_{L^2}).
\end{array}
\end{equation}

\vskip 0.25cm

We will use without further comments the elementary observation
\begin{equation}
||\int{f(x-s,k)g(s)ds} ||_{L^2(\mathbb{R}^{2n})} \leqslant ||g||_{L^1(\mathbb{R}^n)} \, ||f||_{L^2(\mathbb{R}^{2n}) }
\end{equation}

\vskip 0.25cm

\begin{lemma}[$W \approx W_1$]\label{thrm007}
 $\forall t \in [0,T]$
\[
 ||W^\e_1(t) -W^\e(t)||_{L^2} =O( \,\, \varepsilon^{\gamma\frac{1+\theta}{2}-1} ||W_0^\varepsilon||_{L^2} \,\,).
\]
\end{lemma}

\noindent {\bf Proof:}  Denote
\begin{equation}
h^\e=W^\e-W^\e_1.
\end{equation}

Obviously,
\begin{equation}
\begin{array}{c}
\partial_t {h} +  2\pi k \cdot \partial_x {h}+ 
T^{\Vtild}h = -T^{V-\Vtild}W^\e_1, \\ { } \\

h(x,k,0)=0.
\end{array}
\end{equation}

Since the Wigner equation has a bounded $L^2$ propagator, it suffices to bound in $L^2$ the rhs. Indeed, we have 
\begin{equation}
\begin{array}{c}
|| T^{V-\Vtild}W^\e_1 ||_{L^2} = C || \int{\widehat{V}(S)(1-e^{-\frac{\pi}{2}\e^{\gamma}S^2})\widehat{W^\e_1}(X-S,K)\frac{sin(\pi \var SK)}{\var}dS}||_{L^2} \leqslant \\ { } \\

\leqslant O(\var^{-1}) || \widehat{V}(S)(1-e^{-\frac{\pi}{2}\e^{\gamma}S^2})||_{L^1} ||W^\e_1||_{L^2} \leqslant \\ { } \\

\leqslant O(\var^{-1}) ||W^\e_1||_{L^2} \left[{ \e^\gamma \int\limits_{|S| \leqslant \e^{-\frac{\gamma}{2}}}{|\widehat{V}(S)| \, |S|^2dS} +  \int\limits_{|S| \geqslant \e^{-\frac{\gamma}{2}}}{|\widehat{V}(S)| dS}}\right]= \\ { } \\

= O(\var^{-1}) ||W^\e_1||_{L^2} \left[{ \e^\gamma+ \e^\gamma \int\limits_{\rho \in (1,\e^{-\frac{\gamma}{2}})}{ \rho^{-\theta}
d\rho} +  \int\limits_{\rho \geqslant \e^{-\frac{\gamma}{2}}}{\rho^{-2-\theta} d\rho}}\right]= \\ { } \\

= O(\var^{-1}) ||W_0^\varepsilon||_{L^2} \left[{\e^\gamma+ \e^{\gamma[-\frac{1}{2}(1-\theta)+1]} +  \e^{\gamma\frac{1+\theta}{2}}}\right].
\end{array}
\end{equation}
Asking that the remainder is small gives
\begin{equation}
\begin{array}{c}
 \gamma>1, \\ { } \\

 \gamma[-\frac{1-\theta}{2}+1]>1 \,\, \Leftrightarrow \,\, \gamma >\frac{2}{1+\theta}, \\ { } \\

 \gamma \frac{1+\theta}{2} >1 \,\, \Leftrightarrow \,\, \gamma >\frac{2}{1+\theta}.
\end{array}
\end{equation}
So finally with the calibration  $\gamma >\frac{2}{1+\theta}$
the proof of Lemma \ref{thrm007} is complete.

\vskip 0.5cm
Moreover:

\begin{lemma}[$W_1 \approx \rho_2^\e$]\label{thrm008}
$\forall t \in [0,T]$
\begin{equation}
 ||W^\e_1(t)-\rho_2^\e(t)||_{L^2} =O(\var^{\frac{\theta}{2+\theta}-\delta} ||W_0^\varepsilon||_{L^2}).
\end{equation}
\end{lemma}

\noindent {\bf Proof:} It is straightforward to check that
\begin{equation}
\begin{array}{c}
\label{eqw1a}
\partial_t \widehat{W}^\e_1 - 2\pi X \cdot \partial_K \widehat{W}^\e_1 +2\int{ \widehat{\Vtild}(S) \widehat{W}^\e_1(X-S,K) \frac{sin(\pi \varepsilon S\cdot K)}{\varepsilon} dS}=0, \\ { } \\

\partial_t\widehat{\rho}^\e_2 - 2\pi X \cdot \partial_K \widehat{\rho}^\e_2 +2\pi\int{ \widehat{\Vtild}(S) \widehat{\rho}^\e_2(X-S,K) S\cdot K dS}=0, \\ { } \\
\end{array}
\end{equation}
and therefore, if $f=\widehat{W^\e}_1-\widehat{\rho}^\e_2$,
\begin{equation}
\begin{array}{c}
\label{eqw1b}
\partial_t f - 2\pi X \cdot \partial_K f +2\int{ \widehat{\Vtild}(S) f(X-S,K) \frac{sin(\pi \varepsilon S\cdot K)}{\varepsilon} dS}=\\
=2\int{ \widehat{\Vtild}(S) \widehat{\rho}^\e_2(X-S,K) \left({ 1-\frac{sin(\pi \varepsilon S\cdot K)}{\pi \varepsilon S \cdot K}}\right)\pi S \cdot K dS}, \\ { } \\
f(x,k,0)=W_0^\varepsilon - [W_0^\varepsilon].
\end{array}
\end{equation}

Using the Duhamel formula,
\begin{equation}
\begin{array}{c}
||f(t)||_{L^2} \leqslant ||W_0^\varepsilon - [W_0^\varepsilon]||_{L^2} + \\ { } \\
+2T \mathop{sup}\limits_{t \in [0,T]} ||\int{ \widehat{\Vtild}(S) \widehat{\rho}^\e_2(X-S,K,t) \left({ 1-\frac{sin(\pi \varepsilon S\cdot K)}{\pi \varepsilon S \cdot K}}\right)\pi S \cdot K dS}||_{L^2}.
\end{array}
\end{equation}
Recall that  the first term above is (relatively) small by assumption; we will work out the other term:

To bound that, first of all observe that
\begin{equation}
 \label{eqBdSinc}
\left|{ \frac{sin(\pi \varepsilon S\cdot K)}{\pi \varepsilon S \cdot K}- 1}\right| \leqslant C min \left\{{ 1,|\pi \varepsilon S \cdot K| }\right\}
\end{equation}
and therefore  for any $b<0$
\begin{equation}
\begin{array}{c}
 \label{eqbdrhs}
|| \int{ \widehat{\Vtild}(S) \widehat{\rho}^\e_2(X-S,K) \left({1- \frac{sin(\pi \varepsilon S\cdot K)}{\pi \varepsilon S \cdot K}}\right)\pi S \cdot K dS }||_{L^2}  \leqslant \\ { } \\

\leqslant C
\var || \int\limits_{0\leqslant |S| \leqslant 1} { |\widehat{\Vtild}(S)| \cdot |S|^2 |\widehat{\rho_2^\e}(X-S,K)|\cdot |K|^2  dS}||_{L^2} + \\
+\var^{1+2b}  ||\int\limits_{1<|S|<\var^b}{\widehat{\Vtild}(S)  \widehat{\rho}^\e_2(X-S,K)|K|^2dS}||_{L^2}+ \\
+||\int\limits_{|S|>\var^b}{\widehat{\Vtild}(S) \widehat{\rho}^\e_2(X-S,K)|SK|dS} ||_{L^2} \leqslant  \\ { } \\

\leqslant C \left({\,\, \var^{1+2b}  ||\widehat{V}(S) min \{ |S|^2,1 \}||_{L^1} ||\rho_2^\e||_{H^2} + \int\limits_{|S|>\var^b}|{\widehat{\Vtild}(S)|\,|S|dS}  ||\rho_2^\e||_{H^1} \,\, }\right) \leqslant \\ { } \\

\leqslant C \left({ \,\, \var^{1+2b}  ||\widehat{V} \frac{S^2}{1+S^2}||_{L^1} ||\rho_2^\e||_{H^2} + \int\limits_{\rho=\var^b}^\infty {\rho^{-1-\theta}d\rho} ||\rho_2^\e||_{H^1} \,\,}\right)=  \\ { } \\

= C \left({ \,\, \var^{1+2b}   ||\rho_2^\e||_{H^2} + \var^{-\theta b} ||\rho_2^\e||_{H^1} \leqslant C   (\var^{1+2b-\delta} + \var^{-\theta b -\delta}) ||W_0^\varepsilon||_{L^2} \,\,}\right)
\end{array}
\end{equation}

Calibrating the parameters is easy; we are given $\theta \in (0,1)$, $\delta$ as in the statement of the Theorem and we need to find $b \in (-\frac{1}2,0)$ so that
\begin{equation}
\begin{array}{c}
1+2b-\delta>0,  \,\,\,\,  -\theta b -\delta >0 \,\,\, \Leftrightarrow \\ 

\Leftrightarrow \frac{\delta-1}{2} < b < -\frac{\delta}\theta.
\end{array}
\end{equation}
It should now be clear that the constraint on $\delta$ in statement of the Theorem comes from the self-consistency check $-\frac{1}2<\frac{\delta-1}{2}  < -\frac{\delta}\theta<0$. The proof of Lemma \ref{thrm008} is complete.

\vskip 0.25cm
Now the last step is to show that $\rho_2^\e \approx \rho_1^\varepsilon$; but this follows by construction, since they satisfy the same equation (which has an $L^2$-continuous propagator), and $||\rho_2^\e(t)-\rho_1^\varepsilon(t)||_{L^2}=||W_0^\varepsilon-[W_0^\varepsilon]||_{L^2}$.

\vskip 0.5cm
The proof of Theorem \ref{DaTheo} is complete.

\subsection{Proof of Theorem \ref{THEO22}} \label{secprtheo22}

One observes of course that, by construction,
\[
|E^\var(\phi)| \leqslant || \rho_1^\varepsilon(t)-W^\varepsilon(t)||_{L^2} ||\phi||_{L^2}. 
\]
This will be useful in situations when we use a test-function from $\mathcal{A} \cap L^2$.

\vskip 0.2cm
Since in general there are  $\phi \in \mathcal{A}\setminus L^2$, we will also  use  a variation of the  space of test-functions $\mathcal{A}$  introduced earlier: let
\begin{equation}
\begin{array}{c}
%\label{eqzatoB}
\mathcal{B}=\{ \phi \in C_0(\mathbb{R}^{2n}) \,|\,\, \int{ \mathop{sup}\limits_{x} |\mathcal{F}_{k \shortrightarrow K}[\phi(x,k)]|dK } + ||\phi||_{L^2} < \infty \}
\end{array}
\end{equation}
equipped with the norm
\begin{equation}
\begin{array}{c}
||\phi||_{\mathcal{B}}=||\phi||_{\mathcal{A}}+||\phi||_{L^2}.
\end{array}
\end{equation}
In particular,  it is clear by construction that $\mathcal{B}\subseteq \mathcal{A}$, $\mathcal{A}'\subseteq \mathcal{B}'$, and $||f||_{L^2} \leqslant ||f||_{\mathcal{B}}$, $||f||_{\mathcal{B}'}\leqslant ||f||_{L^2} $. Moreover, one easily sees the following 

\vskip 0.4cm
\noindent {\bf Claim:} There exists $R>0$ such that $||W^\var(t)-\rho_1^\var(t)||_{\mathcal{A}'}<R$ for all $t\in[0,T]$. 

\vskip 0.2cm
\noindent {\bf Proof of the claim:} Observe that
\[
\begin{array}{c}
|\langle f,\phi \rangle _{\mathcal{A}',\mathcal{A}}| = |\int\limits_{x,K}{ \mathcal{F}_{k\shortrightarrow K} [f(x,k)] \overline{ \mathcal{F}_{k\shortrightarrow K} [\phi(x,k)]} dxdK } | \leqslant  \\ { } \\

\leqslant \int\limits_{K} { \mathop{sup}\limits_{x} |\mathcal{F}_{k\shortrightarrow K} [\phi(x,k)]| \int\limits_{x} |\mathcal{F}_{k\shortrightarrow K} [f(x,k)]| dx dK } \leqslant \\ { } \\

\leqslant 
\int\limits_{K} { \mathop{sup}\limits_{x} |\mathcal{F}_{k\shortrightarrow K} [\phi(x,k)]| dK} \,\,\,

\mathop{sup}\limits_{K} \int\limits_{x} {  |\mathcal{F}_{k\shortrightarrow K} [f(x,k)]| dx},
\end{array}
\]
i.e.
\[
\begin{array}{c}
||f||_{\mathcal{A}'} \leqslant \mathop{sup}\limits_{K} \int\limits_{x}{ \mathcal{F}_{k\shortrightarrow K} [f(x,k)] dx} %\leqslant \\ { } \\

= \mathop{sup}\limits_{K} \int\limits_{x} { | \int\limits_{k}{ e^{-2\pi i kK}f(x,k)dk} | dx} \leqslant ||f||_{L^1}.
\end{array}
\]

To conclude we use the singular value decomposition of the initial data: if $D_0^\var=\sum\limits_{m} \lambda_m |u_m\rangle \langle u_m|$, then
\[
W_0^\var= \sum\limits_m \lambda_m \int{ e^{-2\pi i k y} u^\var_m(x+\frac{\var y}{2})\overline{u}^\var_m(x-\frac{\var y}{2}) dy},
\]
and
\[
W^\var(t)= \sum\limits_m \lambda_m \int{ e^{-2\pi i k y} u^\var_m(x+\frac{\var y}{2},t)\overline{u}^\var_m(x-\frac{\var y}{2},t) dy},
\]
where of course each of $u_m(t)$, $||u_m||_{L^2}=1$, satisfies the Schr\"odinger equation (\ref{eq1}). The point is that the singular values are time-independent and, especially for positive operators, are controlled by the trace:
\[
D_0^\var >0 \, \Rightarrow \, \lambda_m>0 \,\, \forall m \,\,\,\, \Rightarrow \,\,\, tr(D_0^\var)=\sum\limits_{m} \lambda_m =||\lambda_m||_{l^1}.
\]

We readily observe that
\[
\begin{array}{c}
\mathop{sup}\limits_{K} \int\limits_{x}{| \mathcal{F}_{k\shortrightarrow K} [W^\var(x,k)]| dx} =
\mathop{sup}\limits_{K} \int\limits_{x}{ \sum\limits_{m} |\lambda_m u^\var_m(x+\frac{\var K}{2})\overline{u}^\var_m(x-\frac{\var K}{2})| dx } \leqslant tr(D_0^\var),
\end{array}
\]
and of course
\[
||\rho_1^\var(t)||_{L^1} = ||[W^\var_0]||_{L^1}.
\]
The proof of the claim is complete; it is obvious that $R \leqslant tr(D_0^\var)+||[W_0^\var]||_{L^1}$.

\vskip 0.25cm

We will also use the following:

\noindent {\bf Claim:} Let $\chi\in C^\infty(\mathbb{R}^{2n},[0,1])$ be a cutoff function, $\chi(z) =1 \,\, \forall |z|<1$, $\chi(z) =0 \,\, \forall |z|>2$, and
denote $\phi_M(x,k)=\chi(\frac{(x,k)}M)$. Then, 
\[
\begin{array}{c}
||\phi_M||_{L^2} \leqslant  \frac{2^n \pi^{\frac{n}2}}{\Gamma(\frac{n}2+1)} M^{n}   ||\phi||_{\mathcal{A}},
\end{array}
\]
and
\[
\mathop{lim}\limits_{M \shortrightarrow \infty}  ||\phi_M -\phi||_{\mathcal{A}} =0   .
\]

\vskip 0.1cm
\noindent {\bf Proof of the claim:} First of all,  $\chi(x,k) \in \mathcal{S}(\mathbb{R}^{2n}) \subseteq \mathcal{A}(\mathbb{R}^{2n})$.
Moreover, keeping in mind that $||\phi||_{L^\infty} \leqslant ||\phi||_{\mathcal{A}}$,
\[
||\phi_M||_{L^2} \leqslant   ||\phi||_{L^\infty} \sqrt{ \int\limits_{|(x,k)|<2M}{dxdk} } \leqslant  \frac{2^n \pi^{\frac{n}2}}{\sqrt{ \Gamma(n+1) }} M^{n}   ||\phi||_{\mathcal{A}}.
\]

The second part, i.e.
\[
\begin{array}{c}
\mathop{lim}\limits_{M \shortrightarrow \infty}  ||\phi_M-\phi||_{\mathcal{A}} =0
\end{array}
\]
follows e.g. from the density of $\mathcal{S}(\mathbb{R}^{2n})$ in $\mathcal{A}$ \cite{LP}. Additional information on the decay of $\phi$ is needed in order to estimate more precisely the rate of convergence.

The proof of the claim is complete.

\vskip 0.25cm

It is clear now how, given a $\phi \in \mathcal{A}$, we can select a family $\{ \phi_M \} \subseteq \mathcal{B}$ such that $||\phi_M-\phi||_{\mathcal{A}} \shortrightarrow 0$ as $M \shortrightarrow \infty$. Moreover set for brevity 
\[
F(\var)=C(\var^\kappa||W_0^\var||_{L^2}+||[W_0^\varepsilon]-W_0^\varepsilon||_{L^2})
\]
for some large enough $C>0$. Now, recalling our remarks in the beginning of the section and the first claim, it follows that 
\begin{equation}
\label{eq7_7}
\begin{array}{c}
|\langle W^\var-\rho_1^\var, \phi \rangle| \leqslant |\langle W^\var-\rho_1^\var, \phi_M \rangle| + |\langle W^\var-\rho_1^\var, \phi_M -\phi \rangle| \leqslant \\ { } \\

\leqslant F(\var)||\phi_M||_{L^2} +  (1+||W_0^\var||_{L^1}) || \phi_M -\phi||_{\mathcal{A}} \leqslant \\ { } \\

\leqslant F(\var) \frac{2^n \pi^{\frac{n}2}}{\sqrt{ \Gamma(n+1) }} M^n + (1+||W_0^\var||_{L^1}) || \phi_M -\phi||_{\mathcal{A}}
\,\,\,\,\, \forall \var,M . 
\end{array}
\end{equation}

The proof is completed by selecting, for each value of $\var$ a ``small enough'' $M$, so that $M=M(\var) \shortrightarrow \infty$, and
\[
\mathop{lim}\limits_{\var\shortrightarrow 0} F(\var)||\phi_{M(\var)}||_{L^2}=0.
\]
In the statement of the Theorem, the calibration $M= F(\var)^{-\frac{1}{2n}}$ is used; in that case it follows that
$||\phi_{M(\var)}||_{L^2} \leqslant  \frac{\frac{2^n \pi^{\frac{n}2}}{\sqrt{ \Gamma(n+1) }}}{\sqrt{F(\var)}}$. 

The proof is complete.

%
%Moreover, given  $\phi_A \in \mathcal{A}$, we can select a family $\{ \phi_m \}_{m\in \mathbb{N}} \subseteq \mathcal{B}$ such that $||\phi_m-\phi_A||_{\mathcal{A}} \shortrightarrow 0$ as $m \shortrightarrow \infty$. Now set for brevity 
%\[
%F(\var)=C(\var^\kappa||W_0^\var||_{L^2}+||[W_0^\varepsilon]-W_0^\varepsilon||_{L^2})
%\]
%for some large enough $C>0$. Now, recalling our remarks in the beginning of the section and the first claim, it follows that 
%\begin{equation}
%\label{eq7_7}
%\begin{array}{c}
%|\langle W^\var-\rho_1^\var, \phi_A \rangle| \leqslant |\langle W^\var-\rho_1^\var, \phi_m \rangle| + |\langle W^\var-\rho_1^\var, \phi_m -\phi_A \rangle| \leqslant \\ { } \\
%
%\leqslant F(\var)||\phi_m||_{L^2} +  (1+||[W_0^\var]||_{L^1}) || \phi_m -\phi_A||_{\mathcal{A}} %\\ { } \\
%
%\,\,\,\,\, \forall \var,m . 
%\end{array}
%\end{equation}
%
%The proof is completed by selecting, for each value of $\var$ a ``small enough'' $m$, so that $m=m(\var) \shortrightarrow \infty$, and
%\[
%\mathop{lim}\limits_{\var\shortrightarrow 0} F(\var)||\phi_{m(\var)}||_{L^2}=0.
%\]
%
%
%For example, given $\var$, we can select the largest index $m=m(\var)$ so that 
%$||\phi_{m(\var)}||_{L^2} \leqslant \frac{1}{\sqrt{F(\var)}}$. It is clear that as $\var \shortrightarrow 0$, $m(\var) \shortrightarrow \infty$, and therefore equation (\ref{eq7_7}) implies that $|\langle W^\var-\rho_1^\var, \phi_A \rangle| = (1+||[W_0^\var]||_{L^1})o(1) \,\,\, \forall \phi_A \in \mathcal{A}$. This calibration is used in the statement of Theorem \ref{THEO22}.
%

\section{Examples}
\label{secExamples}
This section is devoted in understanding better the non-concentration properties that are in the heart of our main results. There are two versions; we isolate them below as ``Property 1'', a condition to be checked for a given problem that is sufficient for Theorem \ref{DaTheo} to apply to it, and the stronger ``Property 2'', which is sufficient for Theorem \ref{THEO22} and its Corollary \ref{coro1}. Concrete examples of natural problems that satisfy either one or both of them are constructed; it must be noted that it's rather easy to satisfy Property 1 to generic initial data and repulsive singularities. Asking that Property 2 holds is quite restrictive, and it seems that specially prepared initial data are needed for that.

We are only concerned here with examples that give rise to ill-posed classical problems. (Otherwise the existing theory for $C^1$ potentials has no problem to provide a complete description of the semiclassical limit).

\subsection{Examples for Theorem \ref{DaTheo}} \label{subsecproof2}

Here we are concerned with Cauchy problems for the Liouville equation,
\[
\left\{ {
\begin{array}{c}
 \partial_t \rho_1^\varepsilon+2\pi k \cdot \partial_x \rho_1^\varepsilon-\frac{1}{2\pi} \partial_x \Vtild \cdot \partial_k \rho_1^\varepsilon =0,\\ { } \\
\rho_1^\varepsilon(t=0)=W_0^\varepsilon,
\end{array}
}\right.
\]
which have the following

\begin{property} \label{cond2}
There exist $T>0$, $\delta \in (0,\frac{\theta}{2+\theta})$ and $[W_0^\varepsilon] \in H^2(\mathbb{R}^{2n})$ such that 
\begin{equation}\label{approxi223}
||[W_0^\varepsilon]-W_0^\varepsilon||_{L^2} =o(||W_0^\varepsilon||_{L^2})
\end{equation}
and the solution of
\begin{equation}
\label{eqmaikol2222}
\left\{ {
\begin{array}{c}
 \partial_t \rho_2^\e+2\pi k \cdot \partial_x \rho_2^\e-\frac{1}{2\pi} \partial_x \Vtild \cdot \partial_k \rho_2^\e =0,\\ { } \\
\rho_2^\e(t=0)=[W_0^\varepsilon],
\end{array}
}\right.
\end{equation}
satisfies
\begin{equation} \label{eq0polanhsuj}
||\rho_2^\e(t)||_{H^2} =O(\varepsilon^{-\delta}||W_0^\varepsilon||_{L^2}),\ \mbox{ uniformly for } t \in [0,T].
\end{equation}
\end{property}
\vskip 0.25cm

We treat $\kappa,\gamma$ as given parameters; of course when checking the property in the context of Theorem \ref{THEO22} they are controlled by the statement of Theorem \ref{DaTheo}.

\vskip 0.25cm

\begin{lemma} \label{lmscalingsonso}
Let $\theta \in (0,1)$, $\delta\in (0,\frac{\theta}{2+\theta})$, $n= 2$, $V(x)=-|x|^{1+\theta}$ with an appropriate smooth cutoff outside $\{ |x|<2\}$.    Assume moreover that there exists
 a function of compact support $w \in H^2(\mathbb{R}^{4})  \cap L^\infty$,   such that
\begin{equation}
W_0^\varepsilon(x,k)=\delta_x ^{-2} \delta_k^{-2} w(\frac{x-(z_1,z_2)}{\delta _x},\frac{k-(z_3,z_{4})}{\delta_k}),
\end{equation}
with $(z_1,z_2)=-L(z_3,z_4)$ (i.e. the wavepacket is ``shot towards zero'').
$||w||_{H^2}$, $||w||_{L^\infty}$, $|z|$ have to be bounded uniformly in $\varepsilon$, $||w||_{H^2}, ||w||_{L^\infty}, |z|=O(1)$; other than that $w$ and $z$ can be allowed to depend on $\varepsilon$.

The small parameters involved $\delta_x,\delta_k,R =o(1)$ are calibrated as follows:
\[
R = \left({ log\left({\frac{1}{\varepsilon}}\right) }\right)^{-\frac{1}{3}},
\]
 as in equation (\ref{baynena}), and
\begin{equation}
\label{eqswtorr}
\begin{array}{c}
\delta_x^{-1} (R+\delta_k) =o(1), \\ { } \\
 (\delta_x+\delta_k)^{-2}  (\delta_k + R)^{-2} =O( \varepsilon ^{-\frac{\delta}{2}}).
\end{array}
\end{equation}

\vskip 0.25cm
Then Property \ref{cond2} is  satisfied. 
\end{lemma}

\vskip 0.2cm

\noindent {\bf Remarks:} 
\begin{itemize}
\item
In fact it follows from the Lemma that property \ref{cond2} holds for {\em any} $\delta \in (0,\frac{\theta}{2+\theta})$.

\item
The difficulties arise from the non-smootheness at zero, and the details of a smooth behaviour away from zero are irrelevant here. So we will work for a small enough time, before any trajectory starting in a neighbourhood of zero reaches the support of the cutoff function. Therefore, to keep the presentation simple, we will not introduce any explicit treatment of the cutoff function. This approach is followed in the sequel as well.

\item
An example of a scaling satisfying the above constraints is $\delta_x=\sqrt{\delta_k}=\sqrt{R}$.
\end{itemize}

\vskip 0.2cm

\noindent {\bf Proof:}
Assume without loss of generality that $z=(0,-L,0,1)$. 
Consider $\phi:\mathbb{R} \shortrightarrow [0,1]$ to be a $C^{\infty}$ function such that
\begin{equation}
\label{eqdefgfffff}
\begin{array}{c}
\phi(x)=0, \,\,\,\, |x|<\frac{1}{2}, \\
\phi(x)=1, \,\,\,\, |x|>1.
\end{array}
\end{equation}

 The modified initial data $[W_0^\varepsilon]$ is selected as
\begin{equation}
\label{eq6754}
[W_0^\varepsilon](x,k)=W_0^\varepsilon(x,k)\phi\left({ \frac{|x_1|}{R + 2 L \delta_k} }\right).
\end{equation}

The key claim is that, for $T$ small enough but independent from $\e$,
\[
\left({ \mathop{\bigcup}\limits_{t\in [0,T]}  \mbox{supp} \rho_2^\e(t) }\right) \,\,  \cap   \,\, \{ |x|<R \} =\varnothing
\]
% $\rho_2^\e$ never touches $\{ |x|<R \}$. 
 Indeed $[W_0^\varepsilon]$ has two disjointly supported components. The calibration of the parameters and the cut-off is such that, if $[W_0^\varepsilon]$ was propagated by a Liouville equation corresponding to $V(x)=0$, each of the components would stay on either side of the strip $\{ |x_1|<R\}$ (see Figure \ref{Fig222}). The idea is that since there is a potential driving the flow {\em away} from $ x=0 $, the claim follows. (In fact it is slightly more complicated for this case, since the potential depends on $x_2$ as well: one easily gets an $O(1)$ upper bound for the time it takes for the trajectories to reach $x_2=0$, and then checks that in this  longer time, $x_1$ doesn't have the time to change sign. The difference in an $O(1)$ factor which we absorb in $C$ below).

\begin{figure}[htb!]
\centering
\includegraphics[width=65mm,height=70mm]{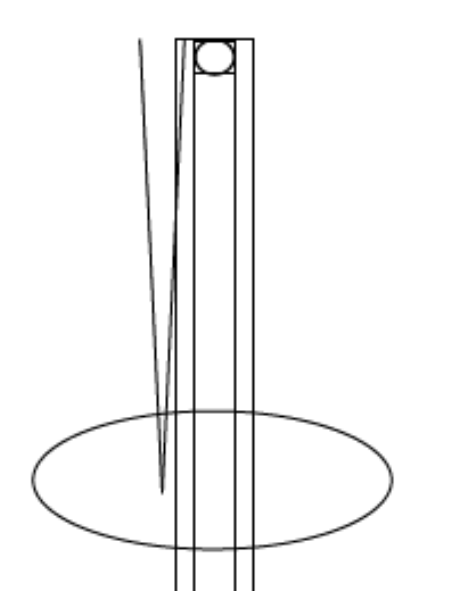}
\caption{ Any trajectory leaving the support of $[W_0^\varepsilon]$ has an initial velocity in a small cone around $k=(0,1)$. If we were in free space, there wouldn't be enough space for sufficient movement in $x_1$ to reach $\{|x_1|<R\}$. Therefore the solution is the sum of two components, supported on either side of $\{|x_1|<R\}=0$. The presence of a repulsive potential in many cases can be easily factored in this construction.} 
\label{Fig222}
\end{figure}

Therefore, in this problem we can take $\VtildEff(x)=\Vtild(x) \phi \left({ \frac{|x|}{R + C L \delta_k} }\right)$. Now for condition (\ref{eq0polanhsuj}), making use of observation \ref{lemJab} (enlarging the cut-off area in fact strengthens observation \ref{lemJab}, i.e. substituting $R + 2 L \delta_k$ for $R$ is painless) and Lemma \ref{lmxf} it follows that
\begin{equation}
\begin{array}{c}
||\rho_2^\e(t)||_{H^2} \leqslant ||[W_0^\varepsilon]||_{H^2} \varepsilon^{-\frac{\delta}{2}} \leqslant C
 ||W_0^\varepsilon||_{H^2} \varepsilon^{-\frac{\delta}{2}} (\delta_k + R)^{-2}  = \\ { } \\
 
= C ||W_0^\varepsilon||_{L^2}  (\delta_x+\delta_k)^{-2}  (\delta_k + R)^{-2} \varepsilon^{-\frac{\delta}{2}},
 \end{array}
\end{equation}
and therefore -- using the scaling of equation (\ref{eqswtorr}) -- $||\rho_2^\e(t)||_{H^2} =O( \varepsilon^{-\delta}||W_0^\varepsilon||_{L^2})$.

Now for $W_0^\varepsilon-[W_0^\varepsilon]$:
\begin{equation}
\begin{array}{c}
||W_0^\varepsilon-[W_0^\varepsilon]||^2_{L^2} \leqslant \\ { } \\

\leqslant \delta_x^{-2n}\delta_k^{-2n} \int\limits_{|x_1|<C(R+\delta_k)}{|w(\frac{(x_1,x_2)}{\delta_x},\frac{(k_1-1,k_2)}{\delta_k})|^2dxdk} = \\ { } \\

=\delta_x^{-n}\delta_k^{-n} \int\limits_{|x_1|<C\frac{R+\delta_k}{\delta_x}}{|w(x_1,x_2,k_1-1,k_2)|^2dxdk} \leqslant \\ { } \\

\leqslant C  ||W_0^\varepsilon||_{L^2}^2 ||w||_{L^\infty}^2 \delta_x^{-1} (R+\delta_k)
\end{array}
\end{equation}

\vskip 0.25cm

The same idea can be applied to different configurations:

\begin{lemma} \label{lmtbswpolmlkbl} Consider the setup of Lemma \ref{lmscalingsonso} with $z=(0,-2,0,1)$ and the only difference that
\begin{equation} \label{razorVeq}
V(x)=-|x_1|^{1+\theta} \psi(x_1)\psi(x_2).
\end{equation}
where $\psi$ is a smooth cutoff function, $\psi=1-\phi \in \mathcal{S}(\mathbb{R})$ ($\phi$ was defined in equation (\ref{eqdefgfffff})).

Then Property \ref{cond2} holds.
\end{lemma}

\vskip 0.25cm
\noindent {\bf Proof:} The geometry is essentially the same as before, and it is clear that the  obvious adaptation of observation \ref{lemJab}, i.e. the one with
\[
\VtildEff(x)=\Vtild(x) \phi(\frac{x_1}{R+C \delta_k})
\]
holds.

\begin{figure}[htb!]
\label{figurezoro}
\centering
\includegraphics[width=55mm,height=85mm]{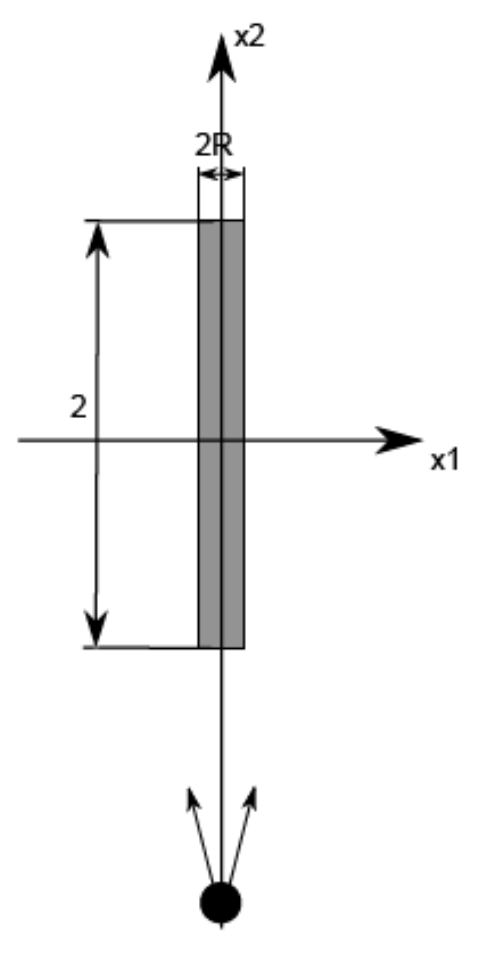}
\caption{ The construction of Lemma \ref{lmscalingsonso} holds, with obvious adjustments, to the example of Lemma \ref{lmtbswpolmlkbl}. } 
\label{Fig9}
\end{figure}

\subsection{Examples for Theorem \ref{THEO22}}
\label{subsecproof4}

In relation to Theorem \ref{THEO22}, we introduce

\begin{property} \label{cond3}
Given $V(x)$, $W_0^e$ with $||W_0^\varepsilon||_{L^2} =o(\varepsilon^{-\kappa})$,
there exist $T>0$, $\delta \in (0,\frac{\theta}{2+\theta})$ and $[W_0^\varepsilon] \in H^2(\mathbb{R}^{2n})$ such that 
\begin{eqnarray}
||[W_0^\varepsilon]-W_0^\varepsilon||_{L^2} =o(1) \label{approxi32}
\end{eqnarray}
and the solution of
\begin{equation}
\label{eqmaikol22221}
\left\{ {
\begin{array}{c}
 \partial_t \rho_2^\e+2\pi k \cdot \partial_x \rho_2^\e-\frac{1}{2\pi} \partial_x \Vtild \cdot \partial_k \rho_2^\e =0,\\ { } \\
\rho_2^\e(t=0)=[W_0^\varepsilon],
\end{array}
}\right.
\end{equation}
satisfies
\begin{equation}
||\rho_2^\e(t)||_{H^2} =O(\varepsilon^{-\delta}||W_0^\varepsilon||_{L^2}),\ \mbox{ uniformly for } t \in [0,T].
\end{equation}
\end{property}
\vskip 0.25cm

As was already seen in section \ref{sspoftlft}, $||f_0^\var-f_0^\var \ast \left({\frac{2}{\var}}\right)^n e^{-2\pi\frac{ x^2+k^2 }{\var}}||_{L^2} =o(1)$. So without loss of generality we can check Property \ref{cond3} in the context of Theorem \ref{thetheo} considering initial data of the form $f_0^\var$, and not the smoothed version that appears there. (More generally, working with compactly supported initial data -- and hence removing the smoothing -- might be very helpful whenever we check Properties \ref{cond2}, \ref{cond3}. In any case, we expect that the smoothing can be removed without loss of generality, just like in this case, whenever the initial data decays rapidly and $||W_0^\var||_{H^1} =o(\var^{-\frac{1}2})$).

\vskip 0.25cm
\begin{lemma} \label{lpbuse9876v} Let $n=1$, $\theta\in(0,1),$ $V(x)=-|x|^{1+\theta}$ with an appropriate smooth cutoff outside $\{ |x|<2 \}$. Assume moreover that $w(x,k) \in H^2 \cap L^\infty \cap L^1$, $supp\, w \subseteq \{ |x|,|k|<1 \}$, and
\begin{equation}
f_0^\varepsilon=\delta_x ^{-1} \delta_k^{-1} w(\frac{x}{\delta _x},\frac{k}{\delta_k}).
\end{equation}
Finally, set $R =  \left({ log\left({\frac{1}{\varepsilon}}\right) }\right)^{-\frac{1}{3}}$.

Then Property \ref{cond3} is satisfied if $\delta_k = C' R^\frac{1+\theta}{2}$, $\delta_x = C'' R^{\frac{1-\theta}{5}}$.
\end{lemma}

\noindent {\bf Proof:} We will cutoff a strip of the form $\{ |x| \leqslant C R \}$, and show that in fact this suffices.

As is illustrated in Figure \ref{fi123oroui}, the preimage under the flow of $\{ |x|<R \,\, \wedge \,\, |k|<\delta_k \}$ is contained between the level sets $\{ 2\pi^2k^2-\Vtild(x)=\Vtild(0)\}$ and $\{ 2\pi^2k^2+\Vtild(x)=\Vtild(R)\}$.% (in red). 

A big help with the algebra will be to approximate $\Vtild(x) =-|x|^{1+\theta} + O(\varepsilon^\frac{\gamma}{2})$. ($\gamma$ is as in the statement of Theorem \ref{DaTheo}). Since all the other small parameters are (negative) powers of $log(\frac{1}{\varepsilon})$, $O(\varepsilon^\frac{\gamma}{2})$ is negligible everywhere. (The justification is that without loss of generality we can localize the problem on a compact set, and then $V(x)\in W^{1,\infty}$ uniformly in $\varepsilon$, see Lemma \ref{lmkjhg}. The conclusion follows by a standard observation on mollifiers, see Lemma \ref{auxlemtkkors}).

It is easily seen that the needed length is $x_*$, defined by
\[
2\pi^2 \delta_k^2+\Vtild(x_*)=\Vtild(R)
\]
or
\[
\begin{array}{c}
|x_*|^{1+\theta} =2\pi^2\delta_k^2+R^{1+\theta}+O(\varepsilon^\frac{\gamma}{2}).
\end{array}
\]

Here is where $\delta_k=C' R^\frac{1+\theta}{2}$ comes from; making that scaling we get $x_*=((1+2\pi^2) R^{1+\theta} + O(\varepsilon^\frac{\gamma}{2}))^{\frac{1}{1+\theta}} = C' R$.

\vskip 0.5cm
So far we have ensured by construction that, if $\delta_k= R^\frac{1+\theta}{2}$, (there is an $O(1)$ constant $C$ so that)  if 
\begin{equation} \label{eqjkjktndh}
[W_0^\varepsilon](x,k)=f_0^\varepsilon(x,k) \phi(\frac{x}{C R}),
\end{equation}
its propagation under equation (\ref{eqmaikol22221}) will never enter $\{ |x|<R \}$. 

Now for the approximation error:
\begin{equation}
\begin{array}{c}
||f_0^\varepsilon-[W_0^\varepsilon]||_{L^2}^2 \leqslant C \delta_x^{-2} \delta_k^{-2} \int\limits_{|x|<CR} {|w(\frac{x}{\delta_x},\frac{k}{\delta_k})|^2dxdk}= \\ { } \\

=C \delta_x^{-1} \delta_k^{-1} \int\limits_{|x|<C\frac{R}{\delta_x}} {|w(x,k)|^2dxdk} \leqslant C
\delta_x^{-2} R^{\frac{1-\theta}{2}}.
\end{array}
\end{equation}
This gives the constraint $R^{\frac{1-\theta}{4}}=o(\delta_x)$ 
for the error to be small, which is satisfied e.g. by our earlier selection $\delta_x = C R^{\frac{1-\theta}{5}}$.

Again, since all the small parameters are powers of $log(\frac{1}{\varepsilon})$, $||[W_0^\varepsilon]||_{H^2}=o(\varepsilon^{-\frac{\delta}{2}})$ follows automatically for any $\delta \in (0,\frac{\theta}{2+\theta})$.

The proof is complete.

\begin{figure}[htb!]
\label{fi123oroui}
\centering
\includegraphics[width=90mm,height=90mm]{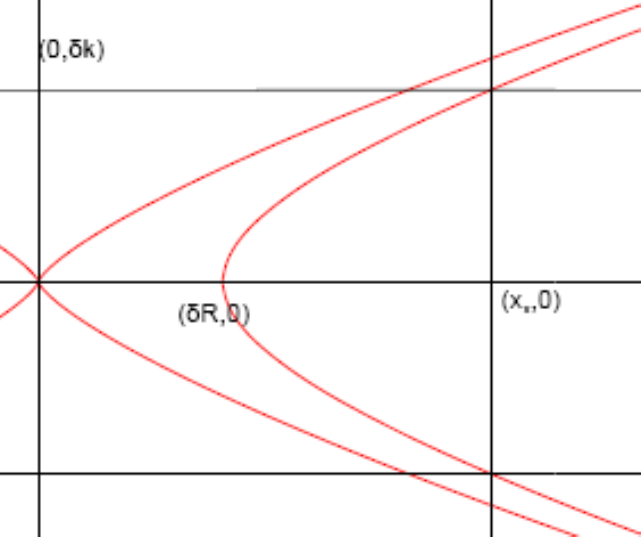}
\caption{ Schematic for Lemma \ref{lpbuse9876v}; cutting off a strip in $x$, $\{ |x|<x_*\}$, suffices to cut off the preimage of $\{ |x|<R \}$} 
\label{Fig10}
\end{figure}

\vskip 0.25cm
The following is a somewhat artificial example, but it highlights an interesting behaviour.  

\begin{lemma} \label{lemthirdartbg}
Assume $n=2$, $V$ as in equation (\ref{razorVeq}). Moreover, there is a function $f \in \mathcal{S}(\mathbb{R}^4)$, $\int{f(x,k)dxdk}=1$, such that
\begin{equation}
supp \,\, f \subseteq \{ |x|<1 \} \times \{ |k|<1  \,\,\wedge \,\, k_1>0 \}.
\end{equation}
Denote by $\bar{f}$ 
\begin{equation}
\bar{f}(x_1,x_2,k_1,k_2)=f(-x_1,x_2,-k_1,k_2)
\end{equation}
Now take
\begin{equation}
\label{eqlabal001}
\begin{array}{c}
W_0^\varepsilon(x,k)=\frac{1}{\delta_x^2 \delta_k^2} 
\left({ c_1 f(\frac{x_1-2R,x_2+2}{\delta_x},\frac{k_1,k_2+1}{\delta_k}) }\right.+ \\ {  } \\

\left.{+c_2 \bar{f}(\frac{x_1-2R,x_2+2}{\delta_x},\frac{k_1,k_2+1}{\delta_k})  }\right) .
\end{array}
\end{equation}

If $\delta_k = C' \delta_x = C'' \sqrt{R}$, this problem has Property \ref{cond3}.
\end{lemma}

\noindent {\bf Proof:}
 Two things should be obvious by construction: that $W_0^\varepsilon \rightharpoonup \delta(x+(0,2),k-(0,1))$, and that the propagation of $W_0^\varepsilon$ under equation (\ref{eqmaikol22}), i.e. $\rho_1^\varepsilon$, is never supported inside $\{ |x_1|<R \}$. In other words, for this specifically constructed data, $W_0^\varepsilon=[W_0^\varepsilon]$.

The scaling of $\delta_x,\delta_k$ in this case controls simply the rate of concentration, and therefore $||W_0^\e||_{H^2}$.
The result follows.

\section{Auxiliary results} \label{secauxres}

\begin{lemma}[$2^{nd}$ order derivatives equations for the Liouville equation] \label{thrmRegLiouv}
Consider the Cauchy problem for the Liouville equation with potential $V(x)$ on $\mathbb{R}^n$,
\begin{equation}
\label{eqlolL1I2O3U34V}
\begin{array}{c}
\partial_t f + 2\pi k \cdot \partial_x f - \frac{1}{2\pi} \partial_x V(x) \cdot \partial_k f=0, \\ { } \\
f(t=0) =f_0.
\end{array}
\end{equation}
There are constants $C_1,C_2>0$ depending only on $n$  such that
\[
||f(t)||_{H^2} \leqslant C_1 e^{ tC_2 \mathop{sup}\limits_{|a|\leqslant 3} ||\partial_x^a V(x)||_{L^\infty} }||f_0||_{H^2}.
\]
\end{lemma}

\noindent {\bf Proof:} The proof follows readily with the method of characteristics.

It suffices to observe that if
\begin{equation}
\label{eqDex2a}
\dot{x}_i(t)=2\pi k_i(t), \,\,\,\,\,\,\,\,\,\,\,\, \dot{k}_i(t)=-\frac{1}{2\pi} \partial_{x_i}{V}(x(t)),
\end{equation}
\begin{equation}
\label{eqDex3b}
\begin{array}{c}
z(t)=f(x(t),k(t)), \\ 
z_{x_i}(t)=\partial_{x_i } f(x(t),k(t)), z_{x_i x_j}(t)=\partial_{x_i x_j} f(x(t),k(t)), \\
z_{k_i }(t)=\partial_{k_i } f(x(t),k(t)), z_{k_i k_j}(t)=\partial_{k_i k_j} f(x(t),k(t)), \\
z_{x_i k_j}(t)=\partial_{x_i k_j} f(x(t),k(t)),
\end{array}
\end{equation}
it follows that
\begin{eqnarray}
\dot{z}=0 \\
\dot{z}_{x_i}(t)=\frac{1}{2\pi} \sum\limits_m \partial_{x_i x_m}{V}(x(t)) \,\, z_{x_i k_m}(t),   \\
\dot{z}_{k_i}(t)=-2\pi z_{x_i}(t), \\
\dot{z}_{x_i x_j}(t)=  \nonumber \\
=\frac{1}{2\pi} \sum\limits_m \left[{  \partial_{x_j x_m}{V} z_{x_i k_m}(t)+ \partial_{x_i x_m}{V} z_{x_j k_m}(t)+ \partial_{x_i x_j x_m}{V} z_{k_m}(t) }\right] \\
\dot{z}_{k_i k_j}(t)=-2\pi (z_{x_j k_i}+z_{x_i k_j}),  \\
\dot{z}_{k_i x_j}(t)=-2\pi z_{x_i x_j}(t)+\frac{1}{2\pi}\sum\limits_m \partial_{x_m x_j}{V}(x(t))z_{k_m k_i}(t).
\end{eqnarray}
 
The result now follows by the Gronwall inequality.

\vskip 1cm

\begin{lemma}\label{auxlemtkkors} Consider a function $f \in W^{1,\infty}$, i.e. $\mathop{sup}\limits_{|a|\leqslant 1} ||\partial_x^a f||_{L^\infty} < \infty$.

Then, if $\widetilde{f}=\left({ \frac{2}{\eta} }\right)^{\frac{n}{2}} \int{ e^{-2\pi \frac{|x-x'|^2}{\eta} }f(x' )dx' }$, we have
\begin{equation}
||f-\widetilde{f}||_{L^\infty} =O(||f||_{W^{1,\infty}}\sqrt{\eta}).
\end{equation}

\end{lemma}

\noindent {\bf Proof:} Take any $\zeta\in (0,\frac{1}{2})$. Now we have
\begin{equation}
\begin{array}{c}
|f(x)-\widetilde{f}(x)|=\left|{\left({\frac{2}{\eta}}\right)^{\frac{n}{2}} \int{ e^{-2\pi \frac{|x'|^2}{\eta}}[f(x-x')-f(x)]dx' } }\right|= \\ { } \\

= \left({\frac{2}{\eta}}\right)^{\frac{n}{2}} \left|{ \int\limits_{|x'|<\eta^{\frac{1}{2}-\zeta}}{ e^{-2\pi \frac{|x'|^2}{\eta}}[f(x-x)'-f(x)]dx' } }\right|+O(\eta^\infty)\leqslant \\ { } \\

\leqslant C ||f||_{W^{1,\infty}} \left({\frac{2}{\eta}}\right)^{\frac{n}{2}} \int\limits_{|x'|<\eta^{\frac{1}{2}-\zeta}}{ e^{-2\pi \frac{|x'|^2}{\eta}}|x'|dx' }+O(||f||_{L^\infty} \eta^\infty) = \\ { } \\
 
 = O(\sqrt{\eta}2^{\frac{n}{2}}) ||f||_{W^{1,\infty}}  \int\limits_{|y|<\eta^{-\zeta}}{ e^{-2\pi |y|^2} |y|dy }+O(||f||_{L^\infty} \eta^\infty) =\\ { } \\
 
 = O(||f||_{W^{1,\infty}}\sqrt{\eta}).
\end{array}
\end{equation}
The proof is complete.

\noindent {\bf Remark:} A sharper version is used in subsection \ref{sspoftlft}, namely
\[
\begin{array}{c}
|f(x)-\widetilde{f}(x)| \leqslant C \left({
\mathop{sup}\limits_{
\begin{scriptsize}\begin{array}{c}
|a|=1 \\
|x-x'|<\eta^{\frac{1}{2}-\zeta}
\end{array}\end{scriptsize}
}
\left|{ \partial_x^af }\right| }\right)\,\, \left({\frac{2}{\eta}}\right)^{\frac{n}{2}} \int\limits_{|x'|<\eta^{\frac{1}{2}-\zeta}}{ e^{-2\pi \frac{|x'|^2}{\eta}}|x'|dx' }+ \\ { } \\

+O(||f||_{L^\infty} \eta^\infty).
\end{array}
\]

\vskip 0.5cm
The following observations are used in Section \ref{secExamples}:

\begin{observation}[Locality of the Liouville equation] \label{ObsHanSolo} For a Liouville equation with initial data of compact support $\rho_0$ and a $C^{1,1}$ potential, interchanging the  potential with any one that coincides with it on the ``path'' of the solution,
\begin{equation}
\label{eqDefS}
\mathbb{S}= \mathop{\bigcup}\limits_{t \in [0,T]} \phi_t( supp \rho_0 ),
\end{equation}
does not change the solution.
\end{observation}

Also,

\begin{lemma}[An $H^2$ estimate] \label{lmxf} Consider equation (\ref{eqmaikol22}) with initial data $\rho_0=[W_0^\varepsilon] \in H^2$. Assume that there is a function $\VtildEff$ such that
\begin{equation}
\Vtild=\VtildEff \,\,\,\,\, on \,\, \mathbb{S},
\end{equation}
where $\mathbb{S}$ is defined as in (\ref{eqDefS}), and
\begin{equation}
\label{eqJaba}
\begin{array}{c}

\mathop{sup}\limits_{|A|\leqslant 3}  || \partial_x^A \VtildEff(x)||_{L^\infty} 
\,\,\,
=O \left({  \,\,\, log \left({ \varepsilon^{-\frac{\delta}{2}} }\right) \,\,\, }\right),
\end{array}
\end{equation}
Then,
\begin{equation}
|| \rho_2^\e ||_{H^2} = O(\,\, \varepsilon^{-\frac{\delta}{2}} ||W_0^\varepsilon||_{H^2}\,\,).
\end{equation}
\end{lemma}

\noindent {\bf Proof:} The proof consists of using $\VtildEff$ in place of $\Vtild$, making use of observation \ref{ObsHanSolo}, and then applying directly Theorem \ref{thrmRegLiouv}.

\vskip 0.25cm
We will also use the following
\begin{lemma} \label{lmkjhg}
It is easy to observe that, if  $\nabla V \in L^\infty$ and $supp W_0^\varepsilon$ is compact, it follows that
\begin{equation}
\exists M>0 :\forall t\in [0,T], \, \varepsilon >0 \,:\, \mathop{\bigcup}\limits_{t\in [0,T]} supp\, \rho_1^\varepsilon(t) \subseteq \{ |(x,k)|<M\}.
\end{equation}
\end{lemma}

\vskip 0.25cm
\begin{observation}[$\VtildEff$ and $R$] \label{lemJab} 
Let $V(x)=-|x|^{1+\theta}$. $\Vtild$ is a mollified version, according to equation (\ref{eqdefVtilde}).

Consider $\phi:\mathbb{R}^n \shortrightarrow [0,1]$ to be a $C^{\infty}$ function such that
\begin{equation}
\label{eqdefgfffff01}
\begin{array}{c}
\phi(x)=0, \,\,\,\, |x|<\frac{1}{2}, \\
\phi(x)=1, \,\,\,\, |x|>1,
\end{array}
\end{equation}
and $|| \phi||_{W^{3,\infty}} \leqslant 10$. It is clear that such a function exists; this is of course an arbitrary requirement, but one that allows us not to carry the cutoff function $\phi$ to other results.

Then, setting $R  \geqslant  \left({ log \left({ \varepsilon^{-1} }\right) }\right)^{-\frac{1}{3}}$, $R=o(1)$, it follows that
\begin{equation}
\label{eq0987}
\VtildEff(x):= \phi(\frac{x}{R})\Vtild(x),
\end{equation}
satisfies
\begin{equation}
\begin{array}{c} \label{imgtkkthfm}
\mathop{sup}\limits_{|A|\leqslant 3}  || \partial_x^A \VtildEff(x)||_{L^\infty} 
\,\,\,
=O( \,\,\, log \left({ \varepsilon^{-\frac{\delta}{2}} }\right) \,\,\,),
\end{array}
\end{equation}
while, of course,
\begin{equation}
\label{eqexp}
{\VtildEff}= {\Vtild} \,\,\, on \,\,\, \mathbb{R}^n \setminus \{ |x|<R  \}.
\end{equation}

\end{observation}

\noindent {\bf Proof:} Observe that $1-\phi(x)=\psi(x) \in \mathcal{S}(\mathbb{R}^n)$. Moreover, making use of  observations \ref{lmkjhg} and \ref{ObsHanSolo}, we can restrict $\Vtild$ to $\{ |x|<M \}$ without loss of generality,
\begin{equation}
\Vtild(x) \, \mapsto \, \Vtild(x)\left({ 1-\phi\left({\frac{x}{2M}}\right) }\right).
\end{equation}
Denote
\begin{equation}
\begin{array}{c}
\VtildEff(x)=\Vtild(x)\left({ 1-\phi\left({\frac{x}{2M}}\right) }\right)\phi\left({\frac{x}{R}}\right)=\\ { } \\
=\Vtild(x)\left({ \psi\left({\frac{x}{2M}}\right)-\psi\left({\frac{x}{R}}\right)  }\right).
\end{array}
\end{equation}
Since
\begin{equation}
|| \partial_x^A \widetilde{f}||_{L^\infty} \leqslant || \partial_x^A {f}||_{L^\infty},
\end{equation}
it suffices to work with
\begin{equation}
\partial_x^A |x|^{1+\theta} \left({ \psi\left({\frac{x}{2M}}\right)-\psi\left({\frac{x}{R}}\right)  }\right).
\end{equation}
Set $x'=\frac{x}{R}$; then $\partial_x=R^{-1} \partial_{x'}$. Now we have
\begin{equation}
\begin{array}{c}
|| \partial_{x}^A |x|^{1+\theta} \left({ \psi\left({\frac{x}{2M}}\right)-\psi\left({\frac{x}{R}}\right)  }\right) ||_{L^\infty} = \\ { } \\

=R^{\theta-2} || \partial_{x}^A |x|^{1+\theta} \left({ \psi\left({\frac{R}{2M}x}\right)-\psi(x)  }\right) ||_{L^\infty} \leqslant \\ { } \\

\leqslant R^{\theta-2} \,\,
 \mathop{sup}\limits_{|A|\leqslant 3} || \partial_x^A x^{1+\theta} ||_{L^\infty([\frac{1}{2},\frac{2M}{R}])} \,\,\,
  \mathop{sup}\limits_{|A|\leqslant 3} || \partial_x^A \left({ \psi\left({\frac{R}{2M}x}\right)-\psi(x)  }\right) ||_{L^\infty(\mathbb{R}^{n})}  \leqslant  \\ { } \\
  
\leqslant (2M)^{1+\theta} ||\phi||_{W^{3,\infty}} R^{-3}.
\end{array}
\end{equation}

Recall that the constant $M$ is chosen so that any trajectory leaving the support of the initial data doesn't not escape $\{ |x|<M \}$ for $t\in [0,T]$. That is, to properly quantify it, one needs to consider initial data of compact support associated with the Liouville equation, and a time-scale $T$. Assuming that the initial data is supported in $\{ |(x,k)|<R_0\}$, it is easy to check that one can set $M = R_0+1+ R_0T+|| \nabla V(x)||_{L^{\infty}}T^2$.
 So now for (\ref{imgtkkthfm}) to hold it suffices that 
\[
\frac{ 2^{1+\theta} \,(1+R_0+R_0T+|| \nabla V(x)||_{L^{\infty}}T^2)^{1+\theta} 10}{C \frac{\delta}2} R^{-3} \leqslant  log(\e^{-1}) 
\]
The constant $C$ comes from the $O(\cdot)$ of equation (\ref{imgtkkthfm}); it is clear that it can be chosen so that the constraint finally becomes
$
 R \geqslant  (log(\e^{-1}) )^{-\frac{1}3},
$
which is satisfied by choosing 
\be 
\label{baynena}
 R =  (log(\e^{-1}) )^{-\frac{1}3}.
 \ee

\vskip 0.2cm
The following is contained in \cite{LP}; we include here a brief mention for completeness:
\begin{lemma}[Density matrices]\label{dktnkmmtl} Let $\mu$ be a probability measure on $\mathbb{R}^{2n}$. Then
\[
W_0^\var(x,k)=\left({\frac{2}\var}\right)^n \int{ e^{-2\pi\frac{(x-x')^2+(k-k')^2}{\var}} d\mu(x',k')}
\]
is the Wigner function of a density matrix, i.e. the corresponding operator $D^\var$ is a positive trace-class operator with $tr(D^\var)=1$.
\end{lemma}

\noindent {\bf Proof:} We know that $tr(D^\var)=\int{W^\var dxdk}=1$. Let us now look at positivity. To that end, observe that
the integral kernel
\[
K^\var(x,y)=\left({\frac{2}\var}\right)^{\frac{n}2} \int\limits_{x_0,k_0}{ e^{2\pi i\frac{k_0}{\var}(x-y)}  e^{-\frac{\pi}{\var}[(x-x_0)^2+(y-x_0)^2]} d\mu(x_0,k_0) }
\]
is the kernel of a positive operator. Indeed:
\[
\begin{array}{c}
\int{K^\var(x,y)u(x)\overline{u}(y)dxdy} = \left({\frac{2}\var}\right)^{\frac{n}2} \int\limits_{x_0,k_0}{ |\langle e^{-2\pi i \frac{k_0}\var x+\frac{\pi}\var (x-x_0)^2} ,u \rangle|^2 d\mu(x_0,k_0) } \geqslant 0.
\end{array}
\]

The proof is complete by observing that the Wigner function corresponding to the kernel $\rho^\var$ is
\[
\int{e^{-2\pi i ky} K^\var(x+\frac{\var y}2,x-\frac{\var y}2)dy}=W_0^\var.
\]

\end{document}